\documentclass{amsart}
\usepackage[all]{xy}
\usepackage{amssymb,amscd}
\usepackage{graphicx,color}

\theoremstyle{plain}

\newtheorem*{thmA}{Theorem A}
\newtheorem*{thmB}{Theorem B}
\newtheorem*{corA}{Corollary A}
\newtheorem*{corB}{Corollary B}
\newtheorem*{lem*}{Lemma}
\newtheorem{theorem}{Theorem}[section]
\newtheorem{lemma}[theorem]{Lemma}
\newtheorem{prop}[theorem]{Proposition}
\newtheorem{cor}[theorem]{Corollary}

\newtheorem*{thm**}{Theorem \ref{thm:dom}}
\newtheorem*{cor**}{Corollary \ref{cor:intro}}
\newtheorem*{thmeq}{Theorem \ref{thm:eq}}

\theoremstyle{remark}
\newtheorem{rem}{Remark}

\newcommand{\adj}{\mathop{{\rm adj}}}

\theoremstyle{definition}

\newtheorem{defn}[theorem]{Definition}

\theoremstyle{remark}

\DeclareMathOperator{\sign}{sign}
\DeclareMathOperator{\Inv}{Inv}
\DeclareMathOperator{\Mult}{Mult}
\DeclareMathOperator{\Th}{Th}
\DeclareMathOperator{\Hom}{Hom}
\DeclareMathOperator{\diam}{diam}
\DeclareMathOperator{\cdim}{cdim}

\DeclareMathOperator{\ncl}{ncl}

\newcommand{\beq}{\begin{equation}}
\newcommand{\eeq}{\end{equation}}
\newcommand{\ov}{\overline}

\newcommand{\BA}{\ensuremath{\mathbb{A}}}

\newcommand{\ben}{\begin{enumerate}}
\newcommand{\een}{\end{enumerate}}

\newcommand{\Tr}{\ensuremath{\mathbb{T}}}

\newcommand{\HH}{\ensuremath{\mathbb{H}}}
\newcommand{\GG}{\ensuremath{\mathbb{G}}}
\newcommand{\D}{\mathcal{D}}

\newcommand{\BN}{\mathbb{N}}
\newcommand{\az}{\mathop{{\ensuremath{\alpha}}}}
\newcommand{\BZ}{\mathbb{Z}}

\newcommand{\factor}[2]{{\raise0.7ex\hbox{$#1$} \!\mathord{\left/ {\vphantom {#1 {#2}}}\right.\kern-\nulldelimiterspace}\!\lower0.7ex\hbox{${#2}$}}}

\title[Algebraic Geometry and the Positive Theory of Partially Commutative Groups]{Elements of Algebraic Geometry and the Positive Theory of Partially Commutative Groups}

\author[M. Casals-Ruiz]{Montserrat Casals-Ruiz}\thanks{The first author is supported by Programa de Formaci\'{o}n de Investigadores del Departamento de Educaci\'{o}n, Universidades e Investigaci\'{o}n del Gobierno Vasco}
\address{Montserrat Casals-Ruiz, Department of Mathematics and Statistics, McGill University, 805 Sherbrooke St. West, Montreal, Quebec H3A 2K6, Canada}

\author[I. Kazachkov]{Ilya V. Kazachkov}\thanks{The second author was supported by Richard.~H.~Tomlinson Fellowship}
\address{Ilya V. Kazachkov, Department of Mathematics and Statistics, McGill University, 805 Sherbrooke St. West, Montreal, Quebec H3A 2K6, Canada}

\begin{document}

\begin{abstract}
The first main result of the paper is  a criterion for a partially commutative group $\GG$ to be a
domain. It allows us to reduce the study of algebraic sets over $\GG$ to the study of irreducible
algebraic sets, and reduce the elementary theory of $\GG$ (of a coordinate group over $\GG$) to
the elementary theories  of the direct factors of $\GG$ (to the elementary theory of coordinate
groups of irreducible algebraic sets).

Then we establish normal forms for quantifier-free formulas over a non-abelian directly
indecomposable partially commutative group $\HH$. Analogously to the case of free groups, we
introduce the notion of a generalised equation and prove that the positive theory of $\HH$ has
quantifier elimination and that arbitrary first-order formulas lift from $\HH$ to $\HH\ast F$,
where $F$ is a free group of finite rank. As a consequence, the positive theory of an arbitrary
partially commutative group is decidable.
\end{abstract}

\maketitle

\section{Introduction}

This paper can be considered as a part of a project the aim of which is to construct algebraic (diophantine) geometry over partially commutative groups, and, more generally, to study the elementary theory of partially commutative groups.

Classical algebraic geometry is concerned with the study of the geometry of sets
of solutions of systems of equations, i.e. the geometry of algebraic sets. Taking the collection
of all algebraic sets as a pre-base of closed sets one gets a topology, known as the Zarsiki
topology.  In the Zariski topology, every closed set is a union (maybe infinite) of algebraic sets. In the case that the ring of coefficients or, equivalently, the Zariski topology is Noetherian, every
closed set $Y$ is a finite union of algebraic sets $Y=Y_1\cup \dots\cup Y_k$. In the case that
$Y_i\nsubseteq Y_j$, $i\ne j$, and $Y_i$ can not be non-trivially presented as a union of
algebraic sets, this decomposition is unique and the sets $Y_1, \dots, Y_k$ are referred to as the
irreducible components of $Y$. In general, however, a finite union of algebraic sets is not necessarily again an algebraic set. In classical algebraic geometry, it suffices to require that the ring of coefficients be a domain. Under these assumptions there exists a one-to-one correspondence between algebraic sets and closed sets. Thus, the study of algebraic sets reduces completely to the study of irreducible algebraic sets.

In \cite{AG1} G.~Baumslag, A.~Miasnikov and V.~Remeslennikov lay down the foundations of algebraic
geometry over groups and introduce group-theoretic counterparts of basic notions from algebraic
geometry over fields. The counterpart to the notion of a Noetherian ring is the notion of an
equationally Noetherian group: a group $G$ is called {\em equationally Noetherian} if every system
$S(X) = 1$ with coefficients from $G$ is equivalent to a finite subsystem $S_0 = 1$, where $S_0
\subset S$, i.e. the algebraic set defined by $S$ coincides with the one defined by $S_0$.
The notion of a domain carries over from rings to groups as follows: a group $G$ is called a
\emph{domain} if for any $x,y \ne 1$ there exists $g\in G$ such that $[x,y^g]\ne 1$.

The notions of equationally Noetherian group and domain, play an analogous role (to their
ring-theoretic counterparts) in algebraic geometry over groups (see \cite{AG1}):
\begin{itemize}
    \item a group $G$ is equationally Noetherian if and only if the Zariski topology is
    Noetherian, in particular every closed set is a finite union of algebraic sets;
    \item if a group $G$ is a domain, then the collection of all algebraic sets is a base for the Zariski
    topology.
\end{itemize}

Our main interest in this paper is algebraic geometry over (free) partially commutative groups.
Partially commutative groups are widely studied in different branches of mathematics and computer science, which explains the variety of names they were given: \emph{graph groups}, \emph{right-angled Artin groups}, \emph{semifree groups}, etc. Without trying to give an account of the literature and results in the field we refer the reader to a recent survey \cite{charney} and the introduction and references in \cite{DKR3}.

Partially commutative groups are linear, see \cite{linear}, hence, equationally Noetherian, see
\cite{AG1}. In \cite{AG1} the authors give several sufficient conditions for a group to be a
domain. In particular, any CSA group is a domain and various group-theoretic constructions preserve
the property of being a domain. However, none of the criteria obtained in \cite{AG1} apply to
the case of partially commutative groups. The major obstacle here is that a partially commutative group may contain a direct product of two free groups.

In Section \ref{sec:dom} we give a criterion for a partially commutative group to be a domain:

\begin{thm**}
Let $\GG$ be a partially commutative group. Then $\GG$ is a domain if and only if
$\GG$ is non-abelian and directly indecomposable.
\end{thm**}

Note that even if a partially commutative group is directly indecomposable, it still may contain a direct product of free groups.

The proof of this theorem is given in Section \ref{sec:domains}. It makes use of the technique of van Kampen diagrams over partially commutative groups, which we present in Section \ref{sec:VK} and the description of centralisers
in partially commutative groups (see Theorem \ref{thm:centr}).

The remaining part of the paper has a model-theoretic flavor. In Section \ref{sec:apag}, using
results from \cite{AG3}, we prove that that the elementary theory of $\GG$ (of a coordinate group
over $\GG$) reduces to the elementary theories  of the direct factors of $\GG$ (to the elementary
theory of coordinate groups of irreducible algebraic sets):
\begin{cor**}
Let $\GG$ be a non-abelian directly indecomposable partially commutative group.
\begin{enumerate}
    \item If $Y=Y_1\cup\dots\cup Y_k$ is an algebraic set over $\GG$, where $Y_1,\dots, Y_k$ are the irreducible components of $Y$, then the elementary theory of the coordinate group $\Gamma(Y)$ of $Y$ is decidable if and only if the elementary theory of $\Gamma(Y_i)$ is decidable for all $i=1,\dots, k$.
    \item If $Y=Y_1\cup\dots \cup Y_k$ and $Z=Z_1\cup \dots \cup Z_l$ are two irreducible algebraic sets,  where $Y_1,\dots, Y_k$  and $Z_1,\dots, Z_l$ are the irreducible components of $Y$ and $Z$, respectively, then $\Gamma(Y)$ is elementary equivalent to $\Gamma(Z)$ if and only if $k=l$ and, after a certain re-enumeration, $\Gamma(Y_i)$ is elementary equivalent to $\Gamma(Z_i)$ for all $i=1,\dots, k$.
\end{enumerate}
\end{cor**}

It is known that coordinate groups of algebraic sets over $\GG$ are separated by $\GG$ (are residually $\GG$), see \cite{AG2}. If a coordinate group $\Gamma$ is a coordinate group of an irreducible set, then $\Gamma$ is discriminated by $\GG$ (is fully residually $\GG$), or equivalently, is universally equivalent to $\GG$. Hence, the class of coordinate groups of irreducible algebraic sets is much narrower and admits a convenient logical description.

In his seminal work \cite{Mak82}, Makanin introduced the notion of a generalised equation. In \cite{Mak84} this notion is used in order to show that the existential theory (the compatibility problem) of free groups and monoids is decidable. Since then this result has been generalised in various ways. In \cite{Schulz} Schulz generalised Makanin's result to the case of systems of equations over a free monoid with regular constraints, and in \cite{ratcons} Diekert, Gutierrez and Hagenah showed the decidability of the compatibility problem for systems of equations over a free group with rational constraints. Using the latter result, Diekert and Lohrey show in \cite{DL} that the existential theory of a certain class of graph products of groups is decidable. Furthermore in \cite{DL2}, the authors show the decidability of the existential theory for an even wider class of groups. A common feature of the results mentioned above is that that they reduce the problem to the one for free groups with rational constraints.

One of the main applications of the decidability of the compatibility problem for free groups is
the decidability of the positive theory of the respective group. In the case of free groups this
is a very well known result. In his paper \cite{Merz}, Merzlyakov performs quantifier elimination for positive formulas over free groups by describing the Skolem functions. Then using the result of
Makanin, \cite{Mak84}, one gets the decidability of the positive theory.

The aim of Sections \ref{sec:mr} and \ref{sec:posth} is to carry over the approach of Merzlyakov and Makanin to the case of partially commutative groups.

In Section \ref{sec:mr}, we prove that any positive quantifier-free formula
over a non-abelian directly indecomposable partially commutative groups is equivalent to a single equation. In order to do so we prove that
\begin{enumerate}
\item for any finite system of equations $S_1(X) = 1, \ldots, S_k(X) = 1$ one can effectively find a single equation $S(X) = 1$ such that the algebraic set defined by the equations $S_1, \ldots, S_k$ and by $S$ coincide for any non-abelian directly indecomposable partially commutative group $\GG$,
\item for any finite set of equations $S_1(X) = 1, \ldots, S_k(X) = 1$ one can effectively find a single equation $S(X) = 1$ such that the union of algebraic sets defined by the equations $S_1, \ldots, S_k$ coincides with the algebraic set defined by $S$ for any non-abelian directly indecomposable partially commutative group $\GG$.
\end{enumerate}
In the case of free groups, the first result is due to Malcev, see \cite{Mal2}, and  in \cite{Mak84} Makanin attributes the second result to Gurevich. These results hold in fact in a much more general setting (for groups that satisfy certain first-order formulas), for example in \cite{IFT} it is proven that this is the case for torsion-free, non-abelian, CSA groups that satisfy the Vaught's conjecture, in particular, for all non-abelian fully residually free groups and torsion-free hyperbolic groups. Note, that a non-abelian directly indecomposable partially commutative group is almost never a CSA group. We generalise the results of Malcev and Gurevich to the case of partially commutative groups. The exposition in this section as well as in Section \ref{sec:posth} is based on \cite{IFT}.
As an immediate consequence of these results we get a normal form for first order formulas over partially commutative groups (in fact, over a much wider class of groups).

In Section \ref{sec:posth}, we use the normal form for Van Kampen diagrams obtained in Lemma \ref{lem:prod} to describe the finite number of all possible cancellation schemes for a given equation. This allows us to introduce the notion of a generalised equation for partially commutative groups. Then we introduce an analogue of the, so called, Merzlyakov words and perform quantifier elimination for positive formulas over non-abelian directly indecomposable partially commutative groups.

\begin{thmeq}
If
$$
\GG \models \forall x_1\exists y_1\ldots\forall x_k\exists y_k (S(X,Y,A)=1),
$$
then there exist words {\rm (}with constants from $\GG${\rm)} $q_1(x_1),\ldots , q_k(x_1,\ldots ,x_k) \in \GG[X],$ such that
$$
\GG[X]\models S(x_1, q_1(x_1),\ldots ,x_k, q_k(x_1,\ldots,x_k,A))=1,
$$
i.e. the equation
$$
S(x_1,y_1, \ldots, x_k,y_k,A) = 1
$$
{\rm(}in variables $Y${\rm)} has a solution in the group $\GG[X]$.
\end{thmeq}
Our approach, therefore, is a natural analog of the classical approach of Merzlyakov and Makanin to the positive theory of free groups and avoids the technically involved language of constraints.

In particular, quantifier elimination gives a reduction of the decidability of the positive theory of non-abelian directly indecomposable partially commutative groups to the decidability of the compatibility problem of an equation, which is known to be decidable, see \cite{DM}.

Finally, in order to prove that the positive theory of any partially commutative group is decidable, we need to study the positive theory of the direct product of groups. In folklore, it is known that if $G=H_1\times \cdots \times H_k$, then the positive theory of $G$ is decidable if the positive theories of $H_1, \dots, H_k$ are decidable. However,  we were unable to find a reference till (when this paper was already written) M.~Lohrey pointed out that in \cite{DL}, the authors give a proof of this result. We present another proof of this fact in the Appendix. The proof is purely model-theoretic and makes use of the ideas of the proof of Theorem \ref{thm:elthdir} which is due to Feferman and Vaught, see \cite{FV}.

\section{Preliminaries} \label{sec:prelim}

\subsection{Partially commutative groups}

We begin with the basic notions of the theory of free partially commutative groups.
Recall that a (free) \emph{partially commutative} group is defined as follows. Let $\Gamma$ be a finite, undirected, simple graph. Let $A = V(\Gamma) = \{a_1, \dots , a_n\}$ be the set of vertices of $\Gamma$ and let $F(A)$ be the free group on $A$. Let
$$
R = \{[a_i, a_j] \in F(A) \mid a_i, a_j \in A\hbox{ and there is an edge of $\Gamma$ joining $a_i$ to $a_j$}\}.
$$
The partially commutative group corresponding to the (commutation) graph $\Gamma$ is the  group $\GG(\Gamma)$ with presentation $\langle A \mid R\rangle$.  This means that the only relations imposed on the generators are commutation of some of the generators. When the underlying graph is clear from the context we write simply $\GG$.

From now on $A=  \left\{ a_1, \ldots , a_r \right\}$ always stands for a finite alphabet, its elements being called \emph{letters}. We reserve the term \emph{occurrence} to denote an occurrence of a letter or of the formal inverse of a letter in a word. In a more formal way, an occurrence is a pair (letter, its placeholder in the word).

For a given word $w$ denote $\az(w)$ the set of letters occurring in $w$. For a word $w \in \GG$, we denote by $ \ov{w}$ a geodesic of $w$. For a  word $w\in \GG$ define $\BA(w)$ to be the subgroup of $\GG$ generated by all letters that do not occur in $\ov w$ and commute with $w$.  The subgroup $\BA(w)$ is well-defined (independent of the choice of a geodesic $\ov w$), see \cite{EKR}. An element $w\in \GG$ is called \emph{cyclically reduced} if the length of $\ov{w^2}$ is twice the length of $\ov{w}$.

For a partially commutative  group $\GG$ consider its non-commutation graph $\Delta$. The vertex set $V$ of $\Delta$ is a set of generators $A$ of $\GG$. There is an edge connecting $a_i$ and $a_j$ if and only if $\left[a_i, a_j \right] \ne 1$. The graph $\Delta$ is a union of its connected components $I_1, \ldots , I_k$. Then
\beq \label{eq:decomp}
\GG= \GG(I_1) \times \cdots \times \GG(I_k).
\eeq

Consider $w \in \GG$ and the set $\az(w)$. For this set, just as above, consider the graph $\Delta (\az(w))$ (it is a full subgraph of $\Delta$ with vertices $\az(w)$). This graph can be either connected or not. If it is connected we will call $w$ a \emph{block}. If $\Delta(\az(w))$ is not connected, then we can split $w$ into the product of commuting words

\begin{equation} \label{eq:bl}
w= w_{j_1} \cdot w_{j_2} \cdots w_{j_t};\ j_1, \dots, j_t \in J,
\end{equation}
where $|J|$ is the number of connected components of $\Delta(\az(w))$ and the word $w_{j_i}$ involves letters from the $j_i$-th connected component. Clearly, the words $\{w_{j_1}, \dots, w_{j_t}\}$ pairwise commute. Each word $w_{j_i}$, $i \in {1, \dots,t}$ is a block and so we refer to presentation (\ref{eq:bl}) as the block decomposition of $w$.

An element $w\in \GG$ is called a least root (or simply, root) of $v\in \GG$ if there exists an integer $0\ne m\in \BZ$ such that $v=w^m$ and there does not exists $w'\in \GG$ and $m'\in \BZ$, $|m'|>1$, such that $w={w'}^{m'}$. In this case we write $w=\sqrt{v}$. By \cite{DK}, partially commutative groups have least roots, that is the root element of $v$ is defined uniquely.

The following proposition reduces the conjugacy problem for arbitrary elements of a partially commutative group to the one for block elements.

\begin{prop}[Proposition 5.7 of \cite{EKR}] \label{prop:57}\
Let $w= w_{1} \cdot w_{2} \cdots w_{t}$ and $v= v_{1} \cdot v_{2} \cdots v_{s}$ be cyclically reduced elements decomposed into the product of blocks. Then $v$ and $w$ are conjugate if and only if $s= t$ and, after some certain index re-enumeration,  $w_{i}$ is  conjugate to $v_{i},\, i= 1, \dots ,t$.
\end{prop}

\begin{cor} \label{cor:prop:57}
Let $w= w_1^{r_1} \cdot w_2^{r_2} \cdots w_t^{r_t}$ and $v= v_1^{l_1} \cdot v_2^{l_2} \cdots v_s^{l_s}$ be cyclically reduced elements decomposed into the product of blocks, where $w_i$ and $v_j$ are root elements, $l_i,r_j \in \BZ$, $i=1,\dots, t$, $j=1,\dots, s$. Then $w$ and $v$ are conjugate if and only if $s=t$ and, after some certain index re-enumeration, $r_i=l_i$ and $w_i$ is conjugate to $v_i$, $i= 1, \dots, t$.
\end{cor}

The next result describes centralisers of elements in partially commutative groups. As the definition of ``being a domain'' relies on the structure of centralisers, we shall make substantial use of the following theorem.

\begin{theorem}[Centraliser Theorem, Theorem 3.10, \cite{DK}] \label{thm:centr} \
Let $w\in \GG$ be a cyclically reduced word, $w=v_1\dots v_k$ be its block decomposition. Then, the centraliser of $w$ is the following subgroup of $\GG$:
\begin{equation} \notag
C(w)=\langle \sqrt{v_1}\rangle \times \cdots \times \langle \sqrt{v_k} \rangle\times \BA(w).
\end{equation}
\end{theorem}

\begin{cor} \label{cor:centr}
For any $w\in \GG$ the centraliser $C(w)$ of $w$ is an isolated subgroup of $\GG$, i.e. $C(w)=C(\sqrt{w})$.
\end{cor}

\subsection{Algebraic Geometry over Groups} \label{sec:ag}

In this section we recall basic notions of algebraic geometry over groups, see \cite{AG1} for details.

For the purposes of algebraic geometry over a group $G$, one has to consider the category of $G$-groups, i.e. groups which contain a designated subgroup isomorphic to the group $G$. If $H$ and $K$ are $G$-groups then a homomorphism $\varphi: H \rightarrow K$ is a $G$-homomorphism if $\varphi(g)= g$ for every $g \in G$. In the category of $G$-groups morphisms are $G$-homomorphisms; subgroups are $G$-subgroups, etc.

Let $G$ be a group generated by a finite set $A$, $F(X)$ be a free group with basis $X = \{x_1, x_2, \ldots  x_n\}$, $G[X] = G \ast F(X)$ be the free product of $G$ and $F(X)$. A subset $S \subset G[X]$  is called {\em a system of equations} over $G$. As an element of the free product, the left side of every equation in $S = 1$ can be written as a product of some elements from $X \cup X^{-1}$ (which are called {\em variables}) and some elements from $A\subset G$ ({\em constants}).

A {\em solution} of the system $S(X) = 1$ over a $G$-group $H$ is a tuple of elements $h_1, \ldots, h_n \in H$ such that every equation from $S$ vanishes at $(h_1, \ldots, h_n)$, i.e. $S(h_1, \ldots, h_n)=1$ in $H$. Equivalently, a solution of the system $S = 1$ over $H$ is a $G$-homomorphism $\phi : G[X] \longrightarrow H$ such that $S\subseteq \ker(\phi)$. Denote by $\ncl(S)$ the normal closure of $S$ in $G[X]$, and by $G_S$ the quotient group $\factor{G[X]}{\ncl(S)}$. Then every solution of $S(X) = 1$ in $H$ gives rise to a $G$-homomorphism $G_S \rightarrow H$, and vice versa. By $V_H(S)$ we denote the set of all solutions in $H$ of the system  $S = 1$ and call it  the {\em algebraic set defined by} $S$.

The normal subgroup of $G[X]$ of the form
$$
R(S) = \{ T(X) \in G[X] \ \mid \ \forall A\in H^n \, (S(A) = 1 \rightarrow T(A) = 1) \}
$$
is called the {\em radical of  the system $S$}. Note that $S\subseteq R(S)$. There exists a one-to-one correspondence between algebraic sets $V_H(S)$ of systems of equations in $G[X]$ and radical subgroups.

The quotient group
$$
G_{R(S)}=\factor{G[X]}{R(S)}
$$
is called the {\em coordinate group} of the algebraic set  $V_H(S)$, and every solution of $S(X) = 1$ in $H$ is a $G$-homomorphism $G_{R(S)} \rightarrow H$.

A $G$-group $H$ is called {\em $G$-equationally Noetherian} if every system $S(X) = 1$ with coefficients from $G$ is equivalent over $G$ to a finite subsystem $S_0 = 1$, where $S_0 \subset S$, i.e. the systems $S$ and $S_0$ define the same algebraic set. If $G$ is $G$-equationally Noetherian, then we say that $G$ is equationally Noetherian. If a $G$-group $H$ is equationally Noetherian every algebraic set $V$ in $G^n$ is a finite union of {\em irreducible components} of $V$.

Let $H$ and $K$ be $G$-groups.  We say that a family of $G$-homomorphisms ${\mathcal F} \subset \Hom_G(H,K)$ {\em $G$-separates} ({\em $G$-discriminates}) $H$ into $K$ if for every non-trivial element $h \in H$ (every finite set of non-trivial elements $H_0 \subset H$) there exists $\phi \in {\mathcal F}$ such that $h^\phi \ne 1$ ($h^\phi \neq 1$ for every $h \in H_0$).
In this case we say that $H$ is {\em $G$-separated} ({\em $G$-discriminated}) into $K$ by $\mathcal{F}$. In the case that $G=1$, we simply say that $H$ is separated (discriminated) by $K$.

A $G$-group $H$ is called a {\em $G$-domain} if for any $x,y \ne 1$ there exists $g\in G$ such that $[x,y^g]\ne 1$. In the case that $G$ is $G$-domain, we say that $G$ is a domain.

\section{Van Kampen Diagrams} \label{sec:VK}
In this section we present some preliminary results on Van Kampen diagrams. We refer the reader to \cite{HS} and \cite{O} for a more detailed account on van Kampen diagrams. Our aim here is to review some basic notions and techniques and apply them to the particular case of partially commutative groups.

\subsection{Van Kampen Diagrams in Partially Commutative Groups}

By van Kampen's Lemma (see \cite{HS}) the word $w$ represents the trivial element in a fixed group $G$ given by the presentation $\langle A \mid R \rangle$ if and only if there exists a finite connected, oriented, based, labeled, planar graph $\mathcal{D}$ where each oriented edge is labeled by a letter in $A^{\pm 1}$, each bounded region (cell) of $\mathbb{R}^2 \setminus \mathcal{D}$ is labeled by a word in $R$ (up to shifting cyclically or taking inverses) and $w$ can be read on the boundary of the unbounded region of $\mathbb{R}^2 \setminus \mathcal{D}$ from the base vertex. Then we say that $\mathcal{D}$ is a \emph{van Kampen diagram} for the boundary word $w$ over the presentation $\langle A \mid R \rangle$. If $w=uv^{-1}=_{\GG} 1$ we say that $\mathcal{D}$ is a \emph{van Kampen diagram} realising the equality $u=v$. In the event that a van Kampen diagram $\mathcal{D}$ realises the equality $w=\ov w$ we say that $\mathcal{D}$ is a \emph{geodesic van Kampen diagram} for $w$.

\begin{figure}[!h]
  \centering
   \includegraphics[keepaspectratio,height=2.2in]{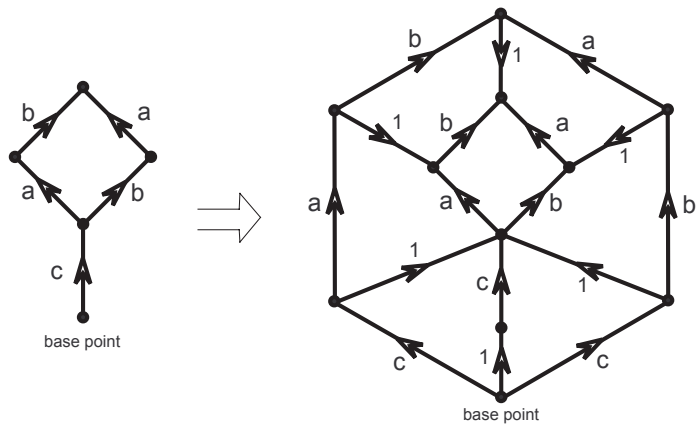}
\caption{van Kampen diagram and non-singular van Kampen diagram for $w= caba^{-1}b^{-1}c^{-1}$ over $\langle a,b,c|[a,b]=1\rangle$.} \label{pic:2}
\end{figure}

Any van Kampen diagram can also be viewed as a 2-complex, with a 2-cell attached for each bounded region (see Figure \ref{pic:2}).

We shall further restrict our considerations to the case when $G$ is a partially commutative group.

Following monograph \cite{O}, if we complete the set of defining relations adding the trivial relations $1\cdot a= a\cdot 1$ for all $a \in A$, then every van Kampen diagram can be transformed so that its boundary is a simple curve. In other words, as a 2-complex the van Kampen diagram is homeomorphic to a disc tiled by cells which are also homeomorphic to a disc (see Figure \ref{pic:2}). We further assume that all van Kampen diagrams are of this form.

Let $\mathcal{D}$ be a van Kampen diagram for the boundary word $w$. Given an occurrence $a$ in $w$, there is a cell $C$ in the 2-complex $\mathcal{D}$ attached to $a$. Since every cell in a van Kampen diagram is either labelled by a relation of the form $a^{-1}b^{-1}ab$ or is a so-called 0-cell, i.e. a cell labelled by $1\cdot a= a\cdot 1$, there is just one occurrence of $a$ and one occurrence of $a^{-1}$ on the boundary of $C$.

Since $\mathcal{D}$ is homeomorphic to a disc, if the occurrence of $a^{-1}$ on the boundary of $C$ is not on the boundary of $\mathcal{D}$, there exists a unique cell $C' \ne C$ attached to this occurrence of $a^{-1}$ in $\mathcal{D}$. Repeating this process, we obtain a unique \emph{band} in $\mathcal{D}$.

\begin{figure}[!h]
  \centering
   \includegraphics[keepaspectratio,height=2in]{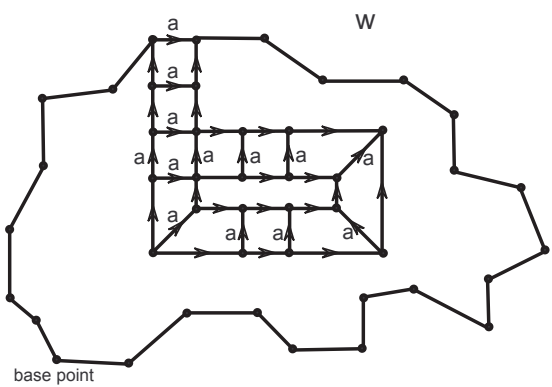}
\caption{Bands do not self-intersect} \label{pic:nointersect}
\end{figure}

Because of the structure of the cells and the fact that $\mathcal{D}$ is homeomorphic to a disc, a band  never self-intersects; indeed, since $\mathcal{D}$ is homeomorphic to a disc, the only way a band can self-intersect is shown in Figure \ref{pic:nointersect}. But then, the cell corresponding to the self-intersection of the band is labelled by the word $aaa^{-1}a^{-1}$.

\begin{figure}[!h]
  \centering
   \includegraphics[keepaspectratio,height=2in]{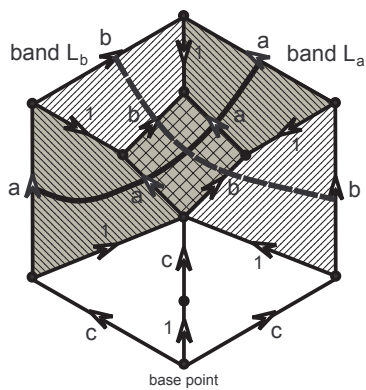}
\caption{Bands in a van Kampen diagram} \label{pic:3}
\end{figure}

Then, since the number of cells in $\mathcal{D}$ is finite, in a finite number of steps the band will again meet the boundary in an occurrence of $a^{-1}$ in $w$ (see Figure \ref{pic:3}).

We will use the notation $L_{a}$ to indicate that a band begins (and thus ends) in an occurrence of a letter $a\in A^{\pm 1}$.

\begin{rem}\label{rem:cross}\

\begin{itemize}
    \item If two bands $L_a$ and $L_b$ cross then the intersection cell realises the equality $a^{-1}b^{-1}ab=1$ and so $a\ne b$ and $[a,b]=1$ {\rm(}see {\rm Figure \ref{pic:3})}.
    \item Every band $L_a$ gives a decomposition of $w$ in the following form $w=w_1aw_2a^{-1}w_3$, where $[a,\az(w_2)]=1$ and $a\notin \az(w_2)$.
\end{itemize}
\end{rem}

\begin{lemma}\label{lem:red}
Let $\GG$ be a partially commutative group. A word $w$ in $\GG$ is not geodesic if and only if $w$ contains a subword $aBa^{-1}$ such that $[a,\alpha(B)]=1$, $a\in A^{\pm 1}$   if and only if there exists a geodesic van Kampen diagram for $w$ that contains a band $L_a$ with both ends in $w$.
\end{lemma}
A proof of this lemma can be found in \cite{Sh}.

\medskip

It is known (see \cite{EKR}) that if a word $w$ represents the trivial element in $\GG$, it can be reduced to the empty word using commutation relations of letters and free cancellation. This reduction process of $w$ to the empty word induces a pairing of occurrences in the word $w$ that cancel. This pairing is independent of the order in which the letters are freely cancelled.

Lemma \ref{lem:red} reflects a consequence of a deeper fact: there exists a one-to-one correspondence between van Kampen diagrams for $w$ and pairings induced by procedures of reductions of $w$ to the empty word.  Indeed, let $\mathcal{D}$ be a van Kampen diagram for the boundary word $w$. Every band $L_a$ gives a decomposition of the form $w=w_{1,a}aw_{2,a}a^{-1}w_{3,a}$. Let $L_a$ be a band such that the length of $w_{2,a}$ is minimal. Hence, every band $L_b$ with an end in an occurrence $b$ in $w_{2,a}$ can not have the other end in an occurrence $b^{-1}$ in $w_{2,a}$. Thus for every occurrence $b$ in $w_{2,a}$ the band $L_b$ crosses the band $L_a$ and hence $[a,\alpha(w_{2,a})]=1$, $a\notin \az(w_{2,a})$. This implies that  $w=w_{1,a}aw_{2,a}a^{-1}w_{3,a}=w_{1,a}w_{2,a}aa^{-1}w_{3,a}$ and thus there exists a process of reduction of $w$ to the empty word in which the occurrence $a$ is cancelled with the occurrence $a^{-1}$. Collapsing the band $L_a$ in $\mathcal{D}$ we get a van Kampen diagram $\mathcal{D'}$ for the boundary word $w'=w_{1,a}w_{2,a}w_{3,a}$, note that the number of cells in $\mathcal{D'}$ is lower than the number of cells in $\mathcal{D}$. The statement follows by induction.

\begin{figure}[!h]
  \centering
   \includegraphics[keepaspectratio,height=2in]{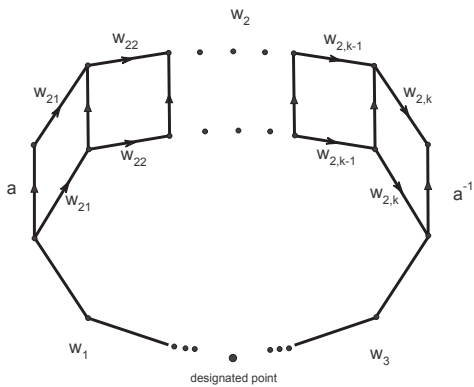}
\caption{Constructing a van Kampen diagram by a process of reduction} \label{pic:11}
\end{figure}

Conversely, if $w$ represents the trivial element in $\GG$, $w$ can be written in the form $w=w_1aw_2a^{-1}w_3$ where $a \in A^{\pm 1}$ and $[a,\alpha(w_2)]=1$ and $a\notin \az(w_2)$. Construct a $|w|$-polygon, designate a point, and orient and label its edges so that starting from the designated point and reading clockwise (or, counterclockwise) one reads $w$. To every edge labelled by an occurrence $w_{2i}$ from $w_2$ we attach a cell labelled by $aw_{2i}a^{-1}w_{2i}^{-1}$. Identifying, as appropriate, the edges labelled by $a^{\pm 1}$ we get a band $L_a$ with ends in $a$ and $a^{-1}$, see Figure \ref{pic:11}. We thereby get a ($|w|-2$)-polygon with the boundary word $w'=w_1w_2w_3$ and thus, by induction, the van Kampen diagram is constructed.

If either in a geodesic van Kampen diagram for $w$ both ends of a band $L_a$ lie in  $w$ or equivalently, if the occurrences $a$ and $a^{-1}$ freely cancel each other in a reduction process of the word $w$ to the empty word, we say that the occurrence $a$ \emph{cancels with} $a^{-1}$. Otherwise, if one of the ends of the band $L_a$ is in an occurrence of $a$ in $w$ and the other is in an occurrence of $a^{-1}$ in $\ov{w}$, we say that $a$ \emph{does not cancel}.

\subsection{Cancellation in a Product of Elements} \label{sec:cancan}

We now consider in detail van Kampen diagrams corresponding to a product of $k$ geodesic words $w_1\cdots w_k=1$.

By Lemma \ref{lem:red} for any van Kampen diagram $\D$ of $w_1\cdots w_k=1$ every band with an end in $w_i$ has its other end in $w_j$, $j\ne i$, $i,j=1,\dots,k$.

\begin{figure}[!h]
  \centering
   \includegraphics[keepaspectratio,width=6in]{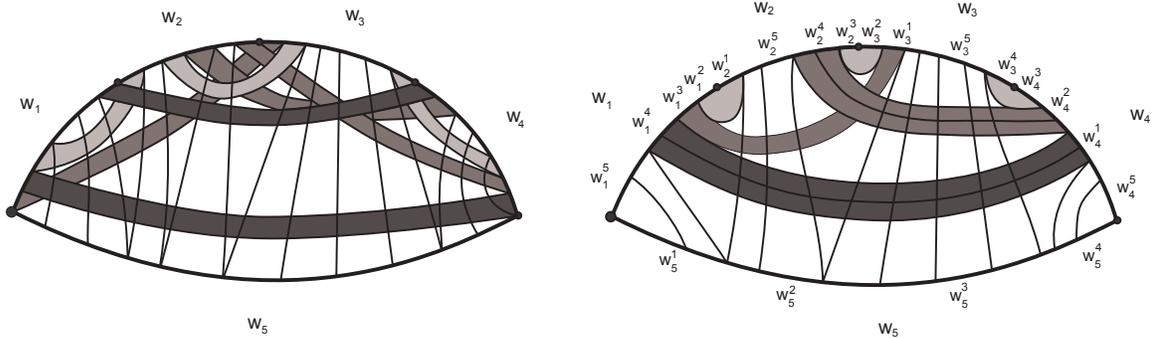}
   \caption{Normal form of a van Kampen diagram, see Lemma \ref{lem:prod}}
\label{pic:8}
\end{figure}

Since every occurrence in $w_1$ cancels, there is a band with an end in a given occurrence $a$ of $w_1$ and another end in $w_i$, $1<i\le k$. Then for any occurrence $b$ in $w_1$ such that
\begin{itemize}
\item $b$ is to the right of $a$, i.e. $w_1=w_1'aw_1'' bw_1'''$ and
\item the band $L_b$ with an end in the occurrence $b$ has its other end in $w_j$, $j>i$,
\end{itemize}
the bands $L_a$ and $L_b$ cross and thus $[a,b]=1$ in $\GG$, see Figure \ref{pic:8}. Therefore the word $w_1$ equals the following geodesic word $w_1=w_1^k\cdots w_1^2$, where the band with an end in any occurrence of $w_1^i$ has its other end in $w_i$.

A similar argument for $w_l$ shows that $w_l$ admits the following decomposition into a  product of geodesic words (perhaps trivial), see Figure \ref{pic:8}:
$$
w_l=w_l^{l-1}\cdots w_l^1 w_l^k\cdots w_l^{l+1},
$$
where the band with an end in any occurrence of $w_l^i$ has its other end in $w_i$.

We summarise the above discussion in the following lemma

\begin{lemma} \label{lem:prod}
Let $\GG$ be a partially commutative group, let $w_1,\dots w_k$ be geodesic words  in $\GG$ such that  $w_1\cdots w_k=1$. Then there exist geodesic words $w_i^j$, $1\le i,j\le k$ such that for any $1\le l\le k$ there exists the following decomposition of $w_l$ into a product of geodesic words:
$$
w_l=w_l^{l-1} \cdots w_l^1 w_l^k\cdots w_l^{l+1},
$$
where $w_l^i={w_i^l}^{-1}$.
\end{lemma}

\begin{cor} \label{cor:prod}
Let $\GG$ be a partially commutative group, let $w_1,\dots w_k,v$ be geodesic words  in $\GG$ such that  $w_1\cdots w_k=v$. Then there exists geodesic words $v_m$, $w_i^j$, $1\le i,j,m\le k$ such that for any $1\le l\le k$ there exists the following decomposition of $w_l$ into a product of geodesic words:
$$
w_l=w_l^{l-1}\cdots w_l^1 v_l w_l^k\cdots w_l^{l+1},
$$
where $w_l^i={w_i^l}^{-1}$ and $v_1\cdots v_k=v$.
\end{cor}

\section{Partially Commutative Groups and Domains} \label{sec:domains}

It is well-known that free groups are domains. The key point of the proof (which relies on the fact that free groups are CSA) is that for $a,x,y\in F$, $x\ne 1$:
$$
\hbox{if } [x,y]=1, [x,y^a]=1,  \hbox{ then } y\in C(a).
$$
Therefore, to see that free groups are domains it suffices to apply the above argument for two elements $a$ and $b$ such that $C(a)\cap C(b)=1$.

Although, directly indecomposable partially commutative groups are not CSA, using the description of centralisers, in Section \ref{sec:dom} we prove that for $a,x,y\in \GG$, such that $x\ne 1$  and $C(a)$ is cyclic:
$$
\hbox{if }  [x,y]=1, [x,y^a]=1,  \hbox{ then either } y\in C(a) \hbox{ or } x\in \BA(y^a).
$$

The aim of Section \ref{sec:canconj} below is to find an element $A\in \GG$ with cyclic centraliser for which $\BA(y^A)=1$.
More precisely, we prove that for any $a \in \GG$, such that $C(a)$ is cyclic, the element $A=a^{2 \cdim(\GG)+2}$ possesses this property. Hence, for $a,x,y\in \GG$, such that $x\ne 1$  and $C(a)$ is cyclic:
$$
\hbox{if } [x,y]=1, [x,y^{A}]=1,  \hbox{ then } y\in C(A).
$$

\subsection{Cancellation and Conjugation}  \label{sec:canconj}

\begin{defn}
We treat the graph $\Delta$ as a metric space with the metric $d$ being the path metric. Let $y$ be a vertex of $\Delta$, define $\adj(y)$ to be $\{v\in \Delta \, \mid \, d(v,y)\le 1 \}$, i.e. the closed ball of radius 1 centered at $y$. For a subset $Y \subseteq A$, set $\adj(Y)=\{v\in \Delta \, \mid \, d(v,y)\le 1 \hbox{ for some } y\in Y \}$.

We set
$$
{\adj}^n(y)=\underbrace{\adj(\adj(\dots\adj(y)\dots))}_{{\small\rm\hbox{$n$ times}}},
$$
thus ${\adj}^n(y)=\{v\in \Delta \, \mid \, d(v,y)\le n \}$ is the closed ball of radius $n$ centered at $y$. Similarly $\adj^n(Y)$, $Y\subseteq A$ is just an $n$-neighbourhood of $Y$ in $\Delta$, $\adj^n(Y)=\{v\in \Delta \, \mid \, d(v,y)\le n \hbox{ for some } y\in Y \}$.

Let $\Delta_1$  be a subgraph of $\Delta$. Then by $\adj(Y)_{\Delta_1}$ we denote the following set $\adj(Y)_{\Delta_1}=\adj(Y)\cap \Delta_1$.
\end{defn}

We shall further use the notion of centraliser dimension $\cdim(G)$ of a group $G$ (see Definition \ref{def:cd} below), an interested reader may consult \cite{DKR1, DKR2, MSh} and references there for a detailed discussion of this notion.

\begin{defn} \label{def:cd}
If there exists an integer $d$ such that the group $G$ has a strictly descending chain of centralisers
$$
C_0>C_1>\dots>C_d
$$
of length $d$ and no centraliser chain of length greater than $d$ then $G$ is said to have {\em centraliser dimension} $\cdim (G)=d$. If no such integer $d$ exists we say that the centraliser dimension of $G$ is infinite, $\cdim (G)= \infty$.
\end{defn}

All partially commutative groups have finite centraliser dimension, \cite{DKR2}.

\begin{lemma} \label{lem:adj}\
Let $\GG$ be a directly indecomposable partially commutative group, $y\in A$ then
\[
{\adj}^{\cdim(\GG)}(y)=A,
\]
i.e. the diameter $\diam(\Delta(\GG))$ of $\Delta(\GG)$ is less than or equal to $\cdim (\GG)$.
\end{lemma}
\begin{proof}
The group $\GG$ is directly indecomposable, hence the non-commutation graph $\Delta(\GG)$ is connected. Therefore for any pair of vertices  $g,h\in V(\Delta)$ there exists a path $p$ of minimal length connecting them. We claim that the length of $p$ is less or   equals $\cdim(\GG)$.

Let $p=(g_0=g, g_1,\dots, g_{r}=h)$. The path $p$ gives rise to a strictly descending chain of centralisers of length $r$:
  \[
    \GG>C(g_0)>C(g_0,g_1)>\dots >C(g_0,\dots,g_{r-2})>C(g_0,\dots,g_{r-1}).
  \]
Indeed, to see that each of the inclusions above is strict we use the minimality of the path $p$. Suppose
$C(g_0,\dots,g_{i-1})=C(g_0,\dots,g_{i})$ for some $1\le i \le r-1$, then since $g_{i+1}\notin C(g_0,\dots,g_{i})$ we also have
$g_{i+1}\notin C(g_0,\dots,g_{i-1})$. So there exists $0\le j\le i-1$, such that $g_j$ does not commute with $g_{i+1}$, thus the
distance between them is $1$. Then $(g_0, g_1, \dots,g_j, g_{i+1}, \dots,g_{r})$ is a shorter path from $g$ to $h$, contradicting the minimality of the path $p$.

As the length $r$ of any strictly descending chain of centralisers is bounded by $\cdim(\GG)$, so is the distance between any two
points in $\Delta$, so ${\adj}^{\cdim(\GG)}(g)=A$.
\end{proof}

\begin{rem}
Note that the equality {\rm $\diam(\Delta(\GG))=\cdim(\GG)$} can be attained. Set $\GG$ to be, for example, the partially commutative group whose non-commutation graph is a path with an odd number of vertices.
\end{rem}

Given two geodesic words $w,v\in \GG$, if the product $wv$ is again geodesic we write $w\circ v$. Let $g \in \GG$ be a geodesic word. We refer to the decomposition $g=g_1 \circ g_2 \circ g_1^{-1}$, where $g_2$ is cyclically reduced as a \emph{cyclic decomposition of g}.

Given a word $g^n$ we write $g^{(i)}$ when we refer to the $i$-th factor $g$ in the product $g^n=g\cdots g$. Similarly, given an occurrence $a$ in $g$ we write $a^{(i)}$, $1 \leq i \leq n$ to indicate that the occurrence $a$ is in $g^{(i)}$.

\begin{lemma}\label{lem:cancellation}
Let $g$, $z \in \GG$ be geodesic and $g=g_1ag_2$, $a\in A^{\pm 1}$. Let $\mathcal{D}$ be a geodesic van Kampen diagram for $gz$.
If the occurrence $a$ does not cancel, neither does any occurrence $b$ in $g_1$ that belongs to $\adj(a)$.
\end{lemma}
\begin{proof}
If $a$ does not cancel, the band $L_{a}$ with an end in $a$ has the other end in
$\ov{gz}$. Then, for any $b \in g_1$ that cancels, corresponding band $L_{b}$ has one end in this occurrence of $b$ in $g^{(1)}$ and the other in an occurrence ${b}^{-1}$ in $z$. Hence the band $L_{b}$ crosses the band $L_{a}$. By Remark \ref{rem:cross}, $b\ne a$ and $[b,a]=1$, so $b \notin \adj(a)$.
\end{proof}

\begin{lemma}\label{lem:cancellation2}
Let $g, z \in \GG$ be geodesic and let $g=g_1ag_2$, $a\in A^{\pm 1}$. Let $\mathcal{D}$ be a geodesic van Kampen diagram for $z^{g}=gzg^{-1}$. If the occurrence of $a$ in $g$ and the corresponding occurrence $a^{-1}$ in $g^{-1}$ do not cancel, then for $b\in \adj(a)$, no occurrence of $b$ in $g_1$ cancels with $b^{-1}$ in $g_1^{-1}$.
\end{lemma}
\begin{proof}
Proof is analogous to the one of Lemma \ref{lem:cancellation}.
\end{proof}

\begin{cor}\label{cor:cancellation1}
Let $g$, $z \in \GG$ be geodesic and let $g$ be a cyclically reduced block. Let $\mathcal{D}$ be a geodesic van Kampen diagram for $gz$.
If there exists an occurrence $a$ in $g$ that does not cancel, then $g^{\cdim(\GG)+1}z=g \circ z'$, i.e. no occurrence in $g^{(1)}$ cancels.
\end{cor}
\begin{proof}
Let $a$ be an occurrence of $g$ that does not cancel in $gz$.
Since $g^{\cdim(\GG)+1}$ is geodesic, the occurrence  $a^{(\cdim(\GG)+1)}$ in $g^{(\cdim(\GG)+1)}$ does not cancel.

Since $g$ is a block, by definition, the graph $\Delta (\az(g))$ is connected, i.e. the subgroup generated by $\az(g)$ is a directly indecomposable partially commutative group. Thus, applying Lemma \ref{lem:adj} to this subgroup and using the fact that the centraliser dimension of $\langle \az(g) \rangle$ is less than or equal to the centraliser dimension of $\GG$ (see \cite{DKR1}) we get that $\adj(a)_{\Delta(\az(g))}^{\cdim(\GG)}= \az(g)$.

Recursively applying Lemma \ref{lem:cancellation}, we get that no occurrence in $g^{(i)}$ that belongs to $\adj_{\Delta(\az(g))}^{((\cdim(\GG)+1)-i)}(a)$ cancels. Therefore, no occurrence from $g^{(1)}$ cancels in $g^{\cdim(\GG)+1}z$.
\end{proof}

\begin{cor}\label{cor:cancellation}
Let $g$, $z\in \GG$ be geodesic and let $g$ be a cyclically reduced block. Let $\mathcal{D}$ be a geodesic van Kampen diagram for $z^{g}=gzg^{-1}$.
If there exists an occurrence $a$ in $g$, such that $a$ and the corresponding occurrence $a^{-1}$ in $g^{-1}$ do not cancel, then $z^{\left(g^{\cdim(\GG)+1}\right)}=g \circ z' \circ g^{-1}$, where $z'=z^{\cdim(\GG)}$.
\end{cor}
\begin{proof}
The proof is analogous to the one of Corollary \ref{cor:cancellation1}
\end{proof}

\begin{defn} \label{defn:cycperm}
Let $z\in \GG$ be a cyclically reduced word and $g, z_1\in \GG$ be so that $z= g^{-1} \circ z_1$ (in the terminology of \cite{EKR}, $g^{-1}$ is called a left-divisor of $z$). We say that the word $g z g^{-1}=z_1g^{-1}$ is a cyclic permutation of $z$.
\end{defn}

Conjugating a cyclically reduced word $z$ one gets a conjugation of a cyclic permutation of $z$. In particular, all letters of $z$ appear in a geodesic $\ov{z^g}$ for any $g \in \GG$. A more precise description is given in the following lemma.

\begin{figure}[!h]
  \centering
   \includegraphics[keepaspectratio,width=6in]{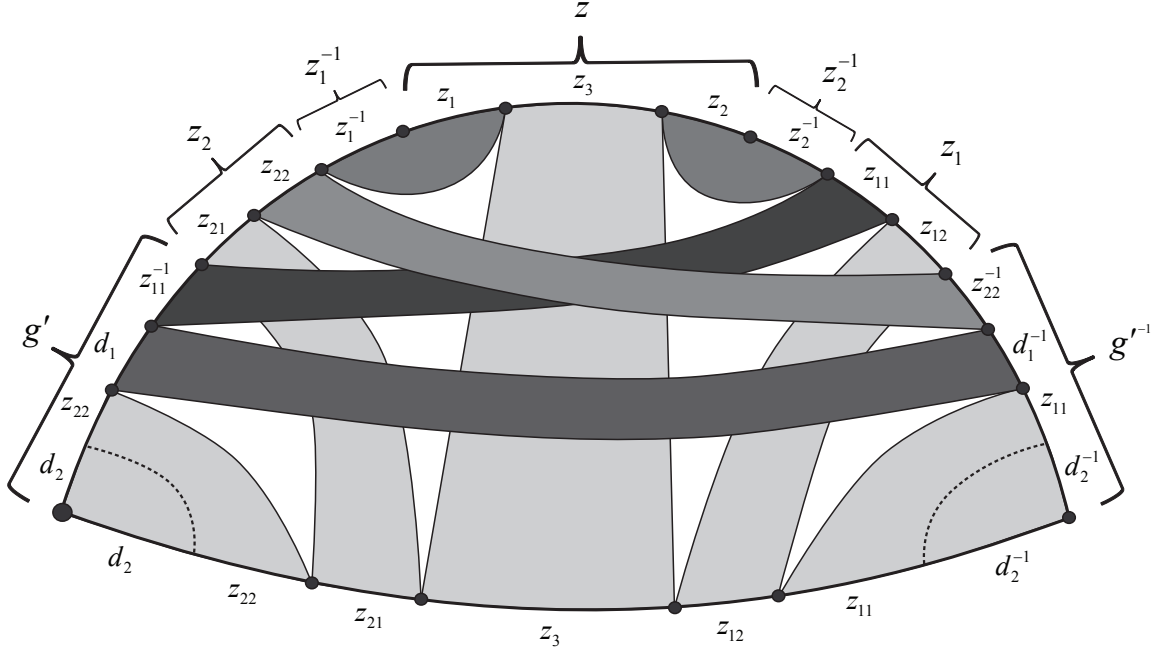}
\caption{Lemma \ref{rem1}: cancellation in $z^g$.} \label{fig:rem1}
\end{figure}

\begin{lemma}\label{rem1}
Let $z,g \in \GG$, let $z$ be cyclically reduced and $g$ be a block. Then there exist decompositions
$$
z=z_1 z_3 z_2,\  z_1=z_{11}z_{12}, \  z_2=z_{21}z_{22}, \ g= d_2d_1z_{22}z_{11}^{-1} z_2 z_1^{-1},
$$
where $z_1\in \BA(z_2)$, $z_{11}\in \BA(z_3)$, $z_{22}\in \BA(z_3)$,
such that $\ov{z^g}= d_2\circ z_{22} \circ z_{21} \circ z_3 \circ z_{12}\circ z_{11} \circ d_2^{-1}$.  Furthermore, if $d_2= 1$ then $d_1=1$.
If $d_2=1$, then either $z_{1}=1$ or $z_2=1$, see {\rm Figure \ref{fig:rem1}}.
\end{lemma}
\begin{proof}
Let $z_1$ be the maximal (where the maximum is taken over all geodesic words representing $g$ and $z$) common initial subword of $z$ and $g^{-1}$, i.e. $z=z_1z_1'$, $g=g_1z_1^{-1}$. Note that such $z_1$ exists and is well-defined, in \cite{EKR} the authors call it the left greatest common divisor of $z$ and $g^{-1}$. Similarly, let $z_2$ be the right greatest common divisor of $z$ and $g$, i.e. $z=z_2'z_2$, $g=g_2z_2$.
Then  $z_1^{-1}$ and $z_2$ are both right divisors of $g$. Let $e$ be the right greatest common divisor of $z_1^{-1}$ and $z_2$. By Proposition 3.18 of \cite{EKR}, $z_1^{-1}=e_1e$, $z_2=e_2e$ and $g=g' e_2{e_1}e$, where $e_1\in \BA(e_2)$. Since $z$ is cyclically reduced, it follows that $e=1$. Hence, $z_1^{-1}=e_1$, $z_2=e_2$, $g=g' z_2{z_1}^{-1}$, $z_1\in \BA(z_2)$.

Apply Corollary \ref{cor:prod} (the notation of which we use below) to the product of words
$$
g'z_2z_1^{-1}zz_2^{-1}z_1{g'}^{-1}=w_1w_2w_3w_4w_5w_6w_7=v=\ov{g'z_2z_1^{-1}zz_2^{-1}z_1{g'}^{-1}}.
$$
Since $g'z_2z_1^{-1}=g$ is geodesic, $w_1^3=w_2^1=w_2^3=w_5^7=w_6^5=w_6^7=1$. By definition of $z_1$ and $z_2$, we get that $w_1^4=w_2^4=w_6^4=w_7^4=1$ and $w_3=w_3^4$, $w_5=w_5^4$, so $w_3^i=w_5^i=1$, $i\ne 4$. As $\az(z_1)\cap \az(z_2)=1$, so $w_2^6=1$.

It follows that $z_2=w_2=v_2w_2^7=z_{21}z_{22}$ and, analogously, $z_1=w_6= w_6^1v_6=z_{11}z_{12}$. Moreover, $g'=w_1=v_1w_1^7w_1^6=v_1w_1^7z_{11}^{-1}$, and ${g'}^{-1}=w_7=w_7^2w_7^1v_7=z_{22}^{-1}w_7^1v_7$.

We obtain that $g'=v_1w_1^7z_{11}^{-1}=v_7^{-1}w_1^7z_{22}$, thus, since $\az(z_{11})\cap \az(z_{22})=1$, we have $v_1=d_2z_{22}$, $v_7^{-1}=d_2z_{11}^{-1}$ and so $g'=d_2z_{22}d_1z_{11}^{-1}=d_2z_{11}^{-1}d_1z_{22}$. This implies that $[d_1,z_{11}z_{22}]=1$, and, since $z_{11} \in \BA(z_{22})$, we have $[d_1,z_{11}]=[d_1,z_{22}]=1$.

By Corollary \ref{cor:prod}, $\ov{z^g}= v_1\cdots v_7=d_2\circ z_{22} \circ z_{21} \circ z_3 \circ z_{12}\circ z_{11} \circ d_2^{-1}$.

Finally, since $g=d_2d_1z_{22}z_{11}^{-1}z_2z_1^{-1}$, if $d_2= 1$ and $d_1\ne 1$, then as shown above $[\az(d_1),\az(z_{11})]=[\az(d_1),\az(z_{22})]=1$, and $[\az(d_1),\az(z_{12})]=[\az(d_1),\az(z_{21})]=1$, see Figure \ref{fig:rem1}. Hence $g$ is not a block. If $d_2=1$, and so, $d_1=1$, then  since $g$ is a block and $[\az(z_1),\az(z_2)]=1$, either $z_{1}=1$ or $z_2=1$.
\end{proof}

\begin{lemma} \label{lem:VK2}
Let $z, g \in \GG$, let $z$ be cyclically reduced, $g$ be a block and let $[z,g]\ne 1$. Furthermore, suppose that $g^{-1}$ does not left-divide $z$ and $z^{-1}$. Then $z^{g^2}=g^2zg^{-2}=g_1\circ z' \circ g_1^{-1}$, where $g_1 \ne 1$ is a left-divisor of $g$, i. e. there exist occurrences $l$ in $g^2$ and, correspondingly, $l^{-1}$ in $g^{-2}$ that do not cancel.
\end{lemma}
\begin{proof}
By and in the notation of Lemma \ref{rem1},  $\ov{z^g}= d_2\circ z_{22} \circ z_{21} \circ z_3 \circ z_{12}\circ z_{11} \circ d_2^{-1}$.

If $d_2\ne 1$, the result follows by Lemma \ref{lem:cancellation2}.

Suppose that $d_2=1$, then, by Lemma \ref{rem1}, $d_1=1$ and either $z_1$ or $z_2$ is trivial. Without loss of generality, we may assume that $z_1=1$. Then $z_2\ne 1$ is a block and
$$
z^{g^2}=g (z_{22}z_{21}z_3)g^{-1}= (z_{22}z_{21}z_{22})( z_{22}z_{21}z_3) (z_{22}z_{21}z_{22})^{-1}.
$$
To prove the lemma it suffices to show that the word $(z_{22}z_{21}z_{22})^{-1}$ does not cancel completely in the product of geodesic words
$$
(z_{22}z_{21}z_{22}z_{22}z_{21}z_3) (z_{22}z_{21}z_{22})^{-1}.
$$

If $\az(z_{21})$ is a subset of $\az(z_3)$, then since $z_{22}\in \BA(z_{3})$, we get that $z_{22} \in \BA(z_{21})$, and $z_2$ is not a block - a contradiction. Therefore, we may write the word $z_{21}^{-1}$ as $A^{-1}B^{-1}$, where  $A^{-1}$ is the left greatest divisor of $z_{21}^{-1}$ such that $\az(A)\subseteq \az(z_3)$ and $B$ is non-trivial. Note that $A\in \BA(z_{22})$ and since $[z_{21},z_{22}]\ne 1$, one has that $[B,z_{22}]\ne 1$.

The word $z^{g^2}$ rewrites as follows:
$$
z^{g^2} = (z_{22}z_{21}z_{22}) z_{22} B  (A z_3 A^{-1}) z_{22}^{-1} B^{-1} z_{22}^{-1}= (z_{22}z_{21}z_{22}) z_{22} B z_3' z_{22}^{-1} B^{-1} z_{22}^{-1}.
$$

By construction, $z_{22}\in \BA(z_{3}')$ and no occurrence in $B^{-1}$ cancels with an occurrence in $z_3'$.

If $[B ,z_{3}']\ne 1$, then there exists an occurrence in $B^{-1}$ that does not cancel and the result follows.
If $[B,z_{3}']= 1$, then we have:
$$
z^{g^2} = ( z_{22}z_{21} z_3'z_{22} z_{22} B) (z_{22}^{-1} B^{-1} z_{22}^{-1}).
$$
If $z_{22}^{-1}$ cancels completely, then we get that both $B$ and $z_{22}$ are right divisors of the word $z_{22}z_{21} z_3'z_{22} z_{22} B$.
Then by Proposition 3.18, \cite{EKR}, $B=B'd$ and $z_{22}=z_{22}'d$, where $z_{22}'\in \BA(B')$. Furthermore, since $[z_{22},B]\ne 1$, then either $[d,z_{22}']\ne 1$ or $[d,B']\ne 1$.

Analogously, if $z_{22}^{-1}B^{-1}$ cancels completely, then we have that both $z_{22}B$ and $Bz_{22}$ are right divisors of the word $z_{22}z_{21} z_3'z_{22} z_{22} B$. Applying again Proposition 3.18, \cite{EKR}, we get that $z_{22}B=U_1D$ and $Bz_{22}=U_2D$ where $U_1 \in \BA(U_2)$.

Combining the above equalities, one has:
$$
z_{22}B=z_{22}'d B'd= U_1 D; \quad
Bz_{22}=B'd z_{22}'d= U_2 D.
$$
Since $\az(U_1) \cap \az(U_2)= \emptyset$, $\az(B') \cap \az(z_{22}')= \emptyset$ and either $[d,z_{22}']\ne 1$ or $[d,B']\ne 1$, the above equalities derive a contradiction. Thus $z_{22}^{-1}B^{-1}$ does not cancel completely in $z^{g^2} = ( z_{22}z_{21} z_3'z_{22} z_{22} B) (z_{22}^{-1} B^{-1} z_{22}^{-1})$ and the result follows.
\end{proof}

We now record some basic properties of $\BA(w)$ which we shall use later. Given $x, y \in \GG$ the following hold:
\begin{enumerate}
    \item[(A)] $x \in \BA(y)$ if and only if $y \in \BA(x)$;
    \item[(B)] if $\az(x) \subset \az(y)$ then $\BA(y)<\BA(x)$;
    \item[(C)] if the centraliser of $x$ is cyclic then $\BA(x)=1$.
\end{enumerate}

\begin{lemma} \label{lem:discr} \
\begin{enumerate}
\item \label{it:lemdis1} Let $g\in \GG$ be a cyclically reduced block and let $z\in \GG$ be so that $g^{-1}$ does not left-divide and right-divide $z$. Then one has $g^{\cdim(\GG)+1}zg^{\cdim(\GG)+1}=g\circ g^{\cdim(\GG)} z g^{\cdim(\GG)} \circ g$.

\item \label{it:lemdis2}Let $g\in \GG$ be a cyclically reduced block and let $z=z_1z_2z_1^{-1}$ be the cyclic decomposition of an element $z\in \GG$. Suppose that $g^{-1}$ does not left-divide $z$, $z^{-1}$, $z_2$ and $z_2^{-1}$, and $[g,z]\ne 1$. Then one has $z^{g^{3\cdim(\GG)+4}}= g \circ z^{g^{3\cdim(\GG)+3}} \circ g^{-1}$.
\end{enumerate}
\end{lemma}
\begin{proof}
We first prove (\ref{it:lemdis1}). We claim that in the product $gzg$ there exist occurrences $l_1$ and $l_2$ in both subwords $g$ of $gzg$ that do not cancel in $\ov{gzg}$. Assume the contrary. Then, one of the subwords $g$ cancels completely. Without loss of generality we may assume that the second subword $g$ of $gzg$ cancels completely. By Lemma \ref{lem:prod}, we get that $g=g_1g_2$, where $g_1$ cancels with $z$ and $g_2$ cancels with $g$. In other words, $z=z'g_1^{-1}$ and $g=g'g_2^{-1}$. Since $g^{-1}$ does not right-divide $z$, we have that $g_2\ne 1$. This derives a contradiction as $g$ is right-divisible by both $g_2$ and $g_2^{-1}$. The statement now follows from Corollary \ref{cor:cancellation1}.

We now prove (\ref{it:lemdis1}). Consider the product $g^{(\cdim(\GG)+1)}z_1$, then no occurrence in $g^{(1)}$ cancels. Indeed, since $g$ does not left-divide $z$, there is an occurrence in $g$ that does not cancel in $gz_1$. Applying Corollary \ref{cor:cancellation1}, we get that no occurrence in $g^{(1)}$ cancels in $g^{(\cdim(\GG)+1)}z_1$.
We thereby get $\overline{g^{(\cdim(\GG)+1)}z_1}= g \circ g' \circ z_1'$, where $z_1'$ is a right-divisor (may be trivial) of $z_1$ and $g'$ is a left-divisor of $g^{\cdim(\GG)}$.

Notice that since $\alpha(g') \subset \alpha(g)$ and $g$ is a block, we get on the one hand that $gg'$ is a block and on the other that for any occurrence $a$ in $g'$ there exists an occurrence $b$ in $g$ such that $b \in \adj(a)$.

If $[z_1', g]\ne 1$ (or $[z_1', g']\ne 1$) then there exists an occurrence $a$ in $g$ (or in $g'$) that belongs to $\adj(\alpha(z_1'))$. Since no occurrence in $z_1'$ and in ${z_1'}^{-1}$ cancels in $z_2^{g g' z_1'}$, by Lemma \ref{lem:cancellation2} neither does the occurrence $a$ in $g$ (or in $g'$), and, correspondingly $a^{-1}$ in $g^{-1}$ (or in ${g'}^{-1}$). If $a$ is an occurrence in $g'$, then there exists an occurrence $b$ in $g$ that belongs to $\adj(a)$, and so by Lemma \ref{lem:cancellation2} this occurrence $b$ in $g$ and the corresponding occurrence $b^{-1}$ in $g^{-1}$ do not cancel. Hence, in any case, there exists an occurrence in $g$ that does not cancel in ${\left({(z_2)}^{ g' z_1'}\right)}^{g}$. Therefore, by Corollary \ref{cor:cancellation} we get that
$$
{\left({(z_2)}^{ g' z_1'}\right)}^{(g^{\cdim(\GG)+1})}=g \circ z_2' \circ g^{-1}
$$
and thus
$$
z^{g^{(2\cdim(\GG)+1)}}={\left({(z_2)}^{g^{(\cdim(\GG)+1)}z_1}\right)}^{(g^{\cdim(\GG)})}=g \circ z^{g^{2\cdim(\GG)}} \circ g^{-1}.
$$

Assume now that $[z_1', g]=1$ and $[z_1',g']=1$. If $[gg',z_2]\ne 1$ or $[g^2g',z_2]\ne 1$, since $gg'$ (or $g^2g'$) is a block, by Lemma \ref{lem:VK2}, there exists an occurrence $a$ in $(gg')^2$ (or in $(g^2g')^2$), and so an occurrence $a$ in $g^{\cdim(\GG)+2}g'$ (or in $g^{\cdim(\GG)+4}g'$), such that $a$ and the occurrence $a^{-1}$ in $(g^{\cdim(\GG)+2}g')^{-1}$ (in ${(g^{\cdim(\GG)+4}g')}^{-1}$) do not cancel in $z_2^{g^{\cdim(\GG)+2}g'}$ (in $z_2^{(g^{\cdim(\GG)+4}g')}$). If $a$ is an occurrence in $g'$, then there exists an occurrence $b$ in $g^{\cdim(\GG)+2}$ (in $g^{\cdim(\GG)+4}$) that belongs to $\adj(a)$, and thus by Lemma \ref{lem:cancellation2} this occurrence $b$ in $g^{\cdim(\GG)+2}$ (in $g^{\cdim(\GG)+4}$) and the corresponding occurrence $b^{-1}$ in $g^{-\cdim(\GG)-2}$ (in $g^{-\cdim(\GG)-4}$) do not cancel.

Now, by Corollary \ref{cor:cancellation} we get
$$
\left({z_2^{g'}}\right)^{g^{(2\cdim(\GG)+4)}}=g \circ z_2' \circ g^{-1}
$$
and thus
\begin{gather} \notag
\begin{split}
z^{g^{(3\cdim(\GG)+4)}}&={\left({(z_2)}^{g^{(2\cdim(\GG)+4)}z_1}\right)}^{(g^{\cdim(\GG)})}={\left({(z_2)}^{g^{\cdim(\GG)+4}g'z_1'}\right)}^{(g^{\cdim(\GG)})}=\\
&={\left({\left({(z_2)}^{g'}\right)}^{(g^{2\cdim(\GG)+4})}\right)}^{z_1'}=g \circ z^{g^{3\cdim(\GG)+3}} \circ g^{-1}.
\end{split}
\end{gather}

Finally, suppose that $[z_1', g]=1$, $[z_1',g']=1$, $[gg',z_2]= 1$ and $[g^2g',z_2]= 1$. We have $g^{\cdim(\GG)+2}=g^2\circ g' \circ d$, $z_1= d^{-1}\circ z_1'$. If $g' \in \langle \sqrt{g} \rangle$, then $d\in \langle \sqrt{g} \rangle$ and so $[d,z_1']=[d,z_2]=1$. This contradicts the assumption that $z=z_2^{d^{-1}z_1'}$ is geodesic. Thus, we may assume that $g' \notin \langle \sqrt{g} \rangle$. In this case, we have that $\az(g')\subseteq \az(g)$, $[g,g']\ne 1$ and $g$ is a block. Therefore, by Theorem \ref{thm:centr}, from $[g,z_1']=[g',z_1']=1$ we get that  $\alpha(g)\subset \BA(z_1')$, and from $[gg',z_2]=[g^2g',z_2]=1$  we get that $\alpha(g)\subset \BA(z_2)$. Since $\alpha(d)\subseteq\alpha(g)$, we have $[d,z_1']=[d,z_2]=1$ -- a contradiction with the assumption that $z=z_2^{d^{-1}z_1'}$ is geodesic.
\end{proof}

\begin{rem}
A more subtle argument shows that the exponent $3\cdim(\GG)+4$ in Lemma \ref{lem:discr} can be replaced by
$\cdim(\GG)$.
\end{rem}

\begin{cor} \label{cor:multconj1}
Let $g\in \GG$ have cyclic centraliser. Then for any element $z\in \GG$ such that $g^{-1}$ does not left-divide $z$, $z^{-1}$, $z_2$ and $z_2^{-1}$, one has  $\BA(z^{(g^{3\cdim(\GG)+4})})=1$.
\end{cor}
\begin{proof}
Let $g=g_1 g_2 g_1^{-1}$ be the cyclic decomposition of $g$. We prove that no occurrence of $g_2^{(1)}$ cancels in $z^{(g^{3\cdim(\GG)+4})}$, therefore, by properties (B) and (C) of $\BA$, we get $\BA(z^{(g^{3\cdim(\GG)+4})})\subseteq \BA(g_2)=1$.

By Lemma \ref{lem:discr}, no occurrence in $g_2^{(1)}$ cancels in ${(z^{{g_1}^{-1}})}^{\left(g_2^{3\cdim(\GG)+4}\right)}$. Then, since
$$
z^{\left(g^{3\cdim(\GG)+4}\right)}={\left({\left(z^{{g_1}^{-1}}\right)}^{g_2^{3\cdim(\GG)+4}}\right)}^{g_1},
$$
no occurrence in $g_2^{(1)}$  cancels.
\end{proof}

\begin{cor} \label{cor:multconj2}
Let $g\in \GG$ have cyclic centraliser. Then for any element $z\in \GG$ such that $\BA(z)\ne 1$ one has $\BA(z^{(g^{3\cdim(\GG)+4})})=1$.
\end{cor}
\begin{proof}
If $g^{-1}$ left-divides $z^{\pm 1}$ or $z_2^{\pm 1}$, by property (B) of $\BA$, we get that $\BA(z)\subset\BA(g)=1$.
\end{proof}

\begin{rem}
In Lemma \ref{lem:discr}, Corollary \ref{cor:multconj1} and Corollary \ref{cor:multconj2}  we impose the condition that $g^{-1}$ does not left-divide $z$, $z^{-1}$, $z_2$ and $z_2^{-1}$, because we seek the bound ${3\cdim(\GG)+4}$ on the number of times one has to conjugate $z$ by $g$. The reason for this is that the notion of centraliser dimension is axiomatisable using (existential) first-order formulas, \cite{DKR1} (we refer the reader to Section \ref{sec:mr} for consequences of this result). If one does not impose this condition, Lemma \ref{lem:discr} could be rephrased as follows.
\end{rem}

\begin{lem*}
Let $g\in \GG$ be a cyclically reduced block, then for any element $z\in \GG$ there exists $N\in \BN$ such that $z^{(g^{N})}=g \circ z^{(g^{N-1})} \circ g^{-1}$.
\end{lem*}

\subsection{Criterion to be a Domain} \label{sec:dom}
\begin{prop} \label{prop:disj}
Let $\GG$ be a non-abelian directly indecomposable partially commutative group. Let $g \in \GG$ have cyclic centraliser, $x, y \in \GG$, $x,y\ne 1$ be such that $[x,y]=1$ and $[x,y^{(g^{3\cdim(\GG)+4})}]=1$. Then $C(x)=C(y)=C(g)$.
\end{prop}

\begin{proof}
Let  $x=w x_1^{r_1}\dots x_k^{r_k} w^{-1}$, where $x_1,\dots,x_k$ are cyclically reduced root elements such that $x_1^{r_1},\dots,x_k^{r_k}$  are the blocks of $x^{w^{-1}}$ and $r_1,\dots,r_k \in \BZ$. Since $[x,y]=1$, by Theorem \ref{thm:centr}, after a certain re-enumeration of indices, we may assume
 \beq \label{eq:y1}
  y=w x_1^{s_1}\dots x_l^{s_l} z w^{-1},
 \eeq where $z\in \BA(x_1\dots x_k)$, $0 \leq l \leq k$ and $s_1,\dots,s_l \in \BZ$. Thus,
 \beq \label{eq:y21}
  y^{(g^{3\cdim(\GG)+4})}={g^{3\cdim(\GG)+4}}w x_1^{s_1}\dots x_l^{s_l} z
  w^{-1}{g^{-(3\cdim(\GG)+4)}}.
 \eeq
Since $[x,y^{(g^{3\cdim(\GG)+4})}]=1$ applying Theorem \ref{thm:centr} once again, we get
  \beq \label{eq:y2}
  y^{(g^{3\cdim(\GG)+4})}=w x_{i_1}^{t_1}\dots x_{i_m}^{t_m} z' w^{-1},
  \eeq
where $z'\in \BA(x_1\dots x_k)$, $0 \leq m \leq k$ and $t_1,\dots,t_m \in \BZ$. Equating (\ref{eq:y21}) and (\ref{eq:y2}), we get
  \beq \label{eq:y3}
  {(x_1^{s_1}\dots x_l^{s_l} z)} ^{w^{-1}{g^{3\cdim(\GG)+4}}w}= x_{i_1}^{t_1}\dots x_{i_m}^{t_m} z'.
  \eeq
Suppose that $l\ge 1$. Then by Corollary \ref{cor:prop:57}, $l=m$ and for any $q \in \{1, \dots, l\}$ there exists $j \in \{1, \dots, m\}$ such that $\az(x_q)=\az(x_{i_j})$, $s_q = t_j$ and $x_q$ is conjugated to $x_{i_j}$ by $w^{-1}{g^{3\cdim(\GG)+4}}w$.

Since $\az(x_q) = \az(x_{i_j})$, and since $x_q$ and $x_{i_j}$ are cyclically reduced root elements whose powers $x_q^{s_q}$ and
$x_{i_j}^{t_j}$ are blocks of the same word $x^{(w^{-1})}$, we get that $x_q=x_{i_j}$ for all $1 \leq q \leq l$, i.e. $y^{w^{-1}}$ and ${(y^{(g^{3\cdim(\GG)+4})})}^{w^{-1}}$ have the same blocks.

From the above it follows that
$$
x_q^{w^{-1}{g^{3\cdim(\GG)+4}}w}=x_{i_j}=x_q,
$$
i. e. $x_q$ commutes with $w^{-1}{g^{3\cdim(\GG)+4}}w$. Since the centraliser of $g$ is cyclic, so is the centraliser of
$w^{-1}{g^{3\cdim(\GG)+4}}w$ and thus, so is the centraliser of $x_q$. More precisely, $C(x_q)={C(g^{3\cdim(\GG)+4})}^{w^{-1}}={C(g)}^{w^{-1}}$.

Since $x_q$ has cyclic centraliser, $x^{w^{-1}}$ and $y^{w^{-1}}$ both have a unique block; furthermore, since $z\in \BA(x_q)$ by property (A) of $\BA$, $z$ is trivial. Therefore, $x=y=(x_q^{r_q})^{w}$ and so $C(x)=C(y)={C(x_q)}^{w}=C(g)$.

Suppose next that $l=0$. We prove then that $x$ is trivial contradicting the assumption. Equations (\ref{eq:y1}) and
(\ref{eq:y3}) rewrite as follows
$$
y=wzw^{-1} \hbox{ and } z'={z}^{w^{-1}{g^{3\cdim(\GG)+4}}w}.
$$
Notice that since $z, z'\in \BA(x_1 \dots x_k)$ by property (C) of $\BA$ we get that $x_1, \dots, x_k \in \BA(z) \cap \BA(z')$.
Therefore, if either $\BA(z)$ or $\BA(z')$ is trivial, so is $x$. Assume $\BA(z)$ is non-trivial. Since the centraliser of $w^{-1}gw$ is cyclic, Corollary \ref{cor:multconj2} applies to $z'={z}^{(({w^{-1}gw})^{{3\cdim(\GG)+4}})}$, thus $\BA(z')=1$ and so
$x=1$.
\end{proof}

\begin{cor} \label{cor:disj}
Let $\GG$ be a non-abelian directly indecomposable partially commutative group. Let $a, b\in \GG$ be elements with cyclic centralisers and such that $C(a)\cap C(b)=1$. Then for any solution $x,y \in \GG$ of the system
$$
[x,y] = 1, \ [x,y^{(a^{{3\cdim(\GG)}+4})}] = 1, \ \ [x,y^{(b^{3\cdim(\GG)+4})}] = 1,
$$
either $x = 1$ or $y = 1$.
\end{cor}
\begin{proof}
Applying Proposition \ref{prop:disj} for the triples $x,y,a$  and $x,y,b$ we get that if $x\ne 1$ and $y \ne 1$, then $C(x)=C(y)=C(a)$ and $C(x)=C(y)=C(b)$ -- a contradiction with $C(a)\cap C(b)=1$. Note, that the elements $a$ and $b$ satisfying the assumption of the corollary exist (it suffices to take two distinct block elements such that $[a,b]\ne 1$ and $\az(a),\az(b)=A$).
\end{proof}

\begin{theorem}[Criterion for a partially commutative group to be a domain]  \label{thm:dom}  \

\noindent
A partially commutative group $\GG$ is a domain if and if $\GG$ is non-abelian and directly indecomposable.
\end{theorem}
\begin{proof}
Since the direct product of two non-abelian groups is never a domain, see \cite{AG1, AG3}, the result follows immediately from Corollary \ref{cor:disj}.
\end{proof}

Note that Corollary \ref{cor:disj} shows in fact that any non-abelian directly indecomposable partially commutative group is a domain with respect to only two elements $a^{3\cdim(\GG)+4}$ and $b^{3\cdim(\GG)+4}$ which are independent of the choice of $x$ and $y$ (in the notation of the definition of domain).

\begin{lemma} \label{lem:for4.17}
Let $a\in \GG$ be a cyclically reduced block element and let $w_1,w_2\in \GG$ be geodesic words of the form
$$
w_1=a^{\delta_ 1} \circ g_1 \circ a^{\epsilon}, \quad w_2=a^{\epsilon} \circ g_2 \circ a^{\delta_2}, \hbox{ where } \epsilon, \delta_1, \delta_2= \pm 1.
$$
Then the geodesic word $\ov{w_1w_2}$  has the form $\ov{w_1w_2}=a^{\delta_1} \circ w_3\circ a^{\delta_2}$.
\end{lemma}
\begin{proof}
We claim that no occurrence from the subword $a^{\epsilon}$ of $w_1$ cancels in the product $w_1w_2$. Indeed, assume the contrary. Let $l$ be the right divisor of $a^{\epsilon}$ of length one that cancels in the product $w_1w_2$, i.e. $a^{\epsilon}=a'l$. It follows that $w_2=l^{-1}w_2'$, see \cite{EKR}. Since $a^{\epsilon}$ and $l^{-1}$ are both left divisors of $w_2$, by Proposition 3.18 of \cite{EKR} one has that either $l^{-1}$ is a left divisor of $a^{\epsilon}$ or  $l^{-1}\in \BA(a)$. As $a^{\epsilon}$ is cyclically reduced and $l$ is a right divisor of $a^{\epsilon}$, so $l^{-1}$ is not a left divisor of $a$. On the other hand, since $l$ is an occurrence of $a^{\epsilon}$, we have that $l^{-1}\notin \BA(a)$. Therefore, no occurrence in the subword $a^{\epsilon}$ of $w_1$ cancels in the product $w_1w_2$.

Since no occurrence in $a^{\epsilon}$ of $w_1$ cancels and $a$ is a block, by Lemma \ref{lem:cancellation}, no occurrence in $a^{\delta_1}$ cancels in the product $w_1w_2$.

An analogous argument shows that no occurrence in the subword $a^{\epsilon}$ of $w_2$ and in $a^{\delta_2}$ cancels and the statement follows.
\end{proof}

\begin{theorem}
Let $\GG$ be a non-abelian directly indecomposable partially commutative group. Then $\GG[X]$ is $\GG$-discriminated by $\GG$.
\end{theorem}
\begin{proof}
The group $\GG[X]$ is a  non-abelian directly indecomposable partially commutative $\GG$-group, thus by Theorem \ref{thm:dom}, $\GG[X]$ is a domain. By Theorem C1 from \cite{BMR2}, it suffices to prove that $\GG[X]$ is $\GG$-separated by $\GG$.

Without loss of generality, we may assume that $X=\{x\}$. Take an element $w\in \GG[X]$. Without loss of generality, we may assume that $w$ is cyclically reduced, $w=x^{k_1}g_1\cdots g_{l-1} x^{k_l}g_{l}$, where $g_i\in \GG$, $g_1, \dots,g_l\ne 1$. Take $a\in \GG$ such that the centraliser of $a$ is cyclic (such $a$ exists since $\GG$ is directly indecomposable) and satisfies the assumptions of Lemma \ref{lem:discr} for every $g_i$, $i=1,\dots, l$, where here, in the notation of Lemma \ref{lem:discr}, $a$ plays the role of $g$ and $g_i$ play the role of $z$.

Consider the homomorphism $\varphi_a:\GG[X]\to \GG$, defined by $x\mapsto a^{6\cdim(\GG)+8}$.  Then
\begin{gather}\notag
\begin{split}
\varphi_a(w)=&a^{k_1(3\cdim(\GG)+4)}\left(a^{k_1(3\cdim(\GG)+4)}g_1a^{k_2(3\cdim(\GG)+4)}\right) \cdot \\ &\left(a^{k_2(3\cdim(\GG)+4)}g_2a^{k_2(3\cdim(\GG)+4)}\right) \cdots \left(a^{k_{l-1}(3\cdim(\GG)+4)}g_{l-1}a^{k_l(3\cdim(\GG)+4)}\right)\cdot a^{k_l(3\cdim(\GG)+4)} g_l
\end{split}
\end{gather}

By Lemma \ref{lem:discr}, every factor of $\varphi_a(w)$ of the form $\left(a^{k_i(3\cdim(\GG)+4)}g_ia^{k_{i+1}(3\cdim(\GG)+4)}\right)$ has the form $a^{\sign(k_i)}\circ \tilde{g_i}\circ a^{\sign(k_{i+1})}$. The statement now follows from Lemma \ref{lem:for4.17}.
\end{proof}

\begin{cor}
Let $\GG$ be a non-abelian directly indecomposable partially commutative group. Then the group $\GG[X]$ is universally equivalent to $\GG$ {\rm(}both in the language of groups and in the language $L_\GG$ enriched by constants from $\GG${\rm)}.
\end{cor}
\begin{proof}
Follows from Theorem C2 in \cite{AG1}.
\end{proof}

\begin{cor}
Let $\GG$ be a non-abelian directly indecomposable partially commutative group. Then
$$
\GG \models \forall X (U(X)=1) \Leftrightarrow \GG[X]\models U(X)=1,
$$
i.e. only the trivial equation has the whole set $\GG^n$ as its solution.
\end{cor}
\begin{proof}
Since $\GG[X]$ is $\GG$-discriminated by $\GG$, if the word $U(X)$ is a non-trivial element of $\GG[X]$, then there exists a $\GG$-homomorphism $\phi:\GG[X]\to\GG$ such that $U^\phi\ne 1$. Then $U(X^\phi)\ne 1$ in $\GG$ -- a contradiction.
\end{proof}

\section{Applications to Algebraic Geometry} \label{sec:apag}

The results and exposition of this section rely on paper \cite{AG3}. We recall here some necessary definitions and restate some results in the case of partially commutative groups. We refer the reader to \cite{AG3} for details and omitted proofs.

A {\em group code} $C$ is a set of formulas
\begin{equation}   \label{eq:code}
    C = \{U(X,P), \ E(X,Y,P),\ \Mult(X,Y,Z,P), \ \Inv(X,Y,P)\}
\end{equation}
where $X, Y, Z, P$ are tuples of variables with $|X| = |Y| = |Z|$. If $P = \emptyset$ then $C$ is called an {\em absolute code} or 0-code.

Let $C$ be a group code, $H$ be a group, and $B$ be an $|P|$-tuple of elements in $H$. We say that $C$ (with parameters $B$) {\em
interprets } a group $C(H,B)$ in $H$ if the following conditions hold:
\ben
 \item [1)] the truth set $U(H,B)$ in $H$ of the formula $U(X,B)$ (with parameters $B$) is non-empty;
 \item [2)] the truth set of the formula $E(X,Y,B)$ (with parameters $B$) defines an equivalence relation $\sim_B$ on $U(H,B)$;
 \item [3)] the formulas $\Mult(X,Y,Z,B)$ and $\Inv(X,Y,B)$ define, correspondingly, a binary  operation ($Z = Z(X,Y)$) and a unary operation ($Y = Y(X)$) on the set $U(H,B)$ compatible with the equivalence relation $\sim_B$;
 \item [4)] the group $C(H,B)$ consists of the set of equivalence classes $\factor{U(H,B)}{\sim_B}$, which form a group with respect to the operations defined by $\Mult(X,Y,Z,B)$ and $\Inv(X,Y,B)$.
\een

We say that a group $G$ is {\em interpretable } (or {\em definable}) in a group $H$ if there exists a group code $C$ and a set of parameters $B \subset H$ such that $G \simeq C(H,B)$. If $C$ is 0-code then $G$ is {\em absolutely} or {\em 0-interpretable} in $H$. The following two types of interpretations are crucial. Let $G$ be a definable subgroup of a group $H$, i.e., there exists a formula $U(x,P)$ and a set of parameters $B \subset H$ such that
$$
 G = \{g \in H \mid H \models U(g,B)\}.
$$
Then $G$ is interpretable in $H$ by the code
$$
 C_G = \{U(x,P), x=y, xy=z, y= x^{-1}\}
$$
with parameters $B$.
If in addition $G$ is a normal subgroup of $H$ then the code
$$
C_{H/G} = \{x=x, \exists v (x = yv \wedge U(v,P)), z= xy, y=x^{-1}\}
$$
interprets the factor-group $\factor{H}{G}$ in $H$ with parameters $B$.
Every group code (\ref{eq:code}) determines a {\em translation} $T_C$ which is a map from the set of all formulas ${\mathcal F}_L$ in the language $L$ into itself. We define $T_C$ by induction as follows:
 \ben
 \item [1)] $T_C(x = y) = E(X,Y,P)$;
 \item [2)] $T_C(xy = z) = \Mult(X,Y,Z,P)$ and $T_C(x^{-1} = y) = \Inv(X,Y,P)$;
 \item [3)] if $\phi, \psi \in {\mathcal F}_L$ and $\circ \in \{ \wedge, \vee, \rightarrow \}$ then
$$
 T_C(\phi \circ \psi) = T_C(\phi) \circ T_c(\psi) \ \ and \ \ T_C(\neg \phi) = \neg T_C(\phi);
$$
 \item [4)] if $\phi \in {\mathcal F}_L$ then
 $$
 T_C(\exists x \phi(x)) = \exists X (U(X,P) \wedge T_C(\phi)),
 $$
 $$
 T_C(\forall x \phi(x)) = \forall X ( U(X,P) \rightarrow T_C(\phi)).
 $$
 \een
Observe, that the formula $T_C(\phi)$ can be constructed effectively from $\phi$.

We say that the elementary theory $\Th(G)$ of a group $G$ is \emph{interpretable} in the group $H$ if there exists a group code $C(H,B)$ of the type (\ref{eq:code}) and a formula $\Psi(P)$ such that $\Th(G) = \Th(C(H,B))$ for any set of parameters $B\subset H$ that satisfies the formula $\Psi(P)$ in $H$.

Any partially commutative group is a direct product of finitely many non-abelian directly indecomposable partially commutative groups and its centre $Z(\GG)$, $Z(\GG)\simeq \BZ^k$, $k\in \BN$. This decomposition is unique up to a permutation of factors. We refer to them as (direct) components of $\GG$.

The centre $Z(\GG)$ is a normal subgroup and a definable subset of $\GG$. It is the truth set of the following formula
\[
\Phi_Z(x): \forall y [x,y]=1,
\]
thus $Z(\GG)$ is 0-interpretable in $\GG$. Consequently, as shown above, the quotient $\factor{\GG}{Z(\GG)}$ is interpretable in $\GG$.

Therefore, to work with partially commutative groups from model-theoretic viewpoint, it suffices to consider free partially
commutative groups with the trivial centre.

Let $\GG$ be a partially commutative group without centre. As mentioned above, in this event $\GG$ is a direct product of directly indecomposable partially-commutative groups, which, in turn are domains by Theorem \ref{thm:dom}. Thus Theorem A and Corollary A of \cite{AG3} apply and can be restated as follows.

\begin{thmA}[cf. \cite{AG3}]
Let $\GG$ be a partially commutative group with trivial centre. Then for each component $\GG_i$ of $\GG$ its elementary theory $\Th(\GG_i)$ is interpretable in the group $\GG$.
\end{thmA}

\begin{corA}[cf. \cite{AG3}]
Let $\GG$ be a partially commutative group and let $\GG=\GG_1\times \dots \times \GG_n\times \BZ^r$, where $\GG_i$ is a non-abelian directly indecomposable partially commutative group, $i=1, \dots, n$. Then the following hold:
\ben
    \item [1)] If $\GG \equiv H$ then
        \begin{itemize}
        \item $H = H_1 \times \ldots \times H_n\times Z(H)$ is a finite direct product of domains and the centre $Z(H)$, with $H_i\equiv \GG_i$ and $Z(H)\equiv Z(\GG)$.
        \item any other decomposition of $H$ as a direct product of domains and its centre has this form {\rm(}after a suitable re-ordering of the factors{\rm)};
            \end{itemize}
    \item [2)] $\Th(\GG)$ is decidable if and only if $\Th(\GG_i)$ is decidable for every $i = 1, \ldots, n$.
\een
\end{corA}

Let $G = G_{1}\times \dots \times G_{k}$ be a direct product of groups $G_i$. A subgroup $H$ of $G$ is called a {\em subdirect product} of groups $G_i$ if $\pi_i(H) = G_i$ for every $i = 1, \ldots, k$, where $\pi_i : G \rightarrow G_i$ is the canonical projection. An embedding
\begin{equation} \label{eq:subdirect1}
\lambda : H \hookrightarrow G_{1}\times \dots \times G_{k}
\end{equation}
is called a {\em subdirect decomposition} of $H$ if  $\lambda(H)$ is a subdirect product of the groups $G_i$. The subdirect decomposition (\ref{eq:subdirect1}) is termed {\em minimal} if $H \cap G_{i}\ne \{ 1 \}$ for every $i = 1, \ldots, k$ (here $G_i$ is viewed as a subgroup of $G$ under the canonical embedding).

\begin{thmB}[\cite{AG3}]
Let $H$ be a minimal subdirect product of domains. Then the elementary theory of each component of $H$ is interpretable in the group $H$.
\end{thmB}

\begin{corB}[\cite{AG3}]
Let $H$ be a minimal subdirect product of $k$ domains and
$$
H \hookrightarrow G_1 \times \ldots \times G_k
$$
be its minimal component decomposition. Then the following hold:
\ben
    \item [1)] if $\Th(H)$ is decidable then $\Th(G_i)$ is  decidable for every $i = 1, \ldots, k$;
    \item [2)] if $\Th(H)$ is $\lambda$-stable then $\Th(G_i)$ is $\lambda$-stable for every $i = 1, \ldots, k$.
\een
\end{corB}

\begin{theorem}\label{th:coord}
Let $\GG$ be a directly indecomposable partially commutative group, and $Y$ be an algebraic set over $\GG$. Then the following conditions hold:
 \ben
 \item [1)] the coordinate group $\Gamma(Y_i)$ of each irreducible component  $Y_i$ of $Y$ is interpretable in the group $\Gamma(Y)$;
 \item [2)] the elementary theory $\Th(\Gamma(Y_i))$ of each irreducible  component $Y_i$ of $Y$ is interpretable in the group $\Gamma(Y).$
 \een
\end{theorem}
\begin{proof} Partially commutative groups are linear (see \cite{H,linear}), thus equationally Noetherian (see \cite{AG1}). We can therefore decompose $Y$ as a finite union of  irreducible algebraic sets, $Y=Y_{1} \cup \dots \cup Y_{k}$, see Corollary 12 in \cite{AG1}. By Proposition 12, \cite{AG1} the coordinate group $\Gamma(Y)$ is a minimal subdirect product of the coordinate groups $\Gamma(Y_1), \dots, \Gamma(Y_k)$.  Every group $\Gamma(Y_i)$, being a coordinate group of an irreducible algebraic set over a domain is again a domain by Theorem D2 in \cite{AG1}.  Now 1), 2) follow from Theorem B.
\end{proof}

\begin{cor}\label{cor:intro}
Let $\GG$ be a directly indecomposable partially commutative group.
\begin{enumerate}
\item If $Y=Y_1\cup\dots\cup Y_k$ is an algebraic set over $\GG$, where $Y_1,\dots, Y_k$ are the irreducible components of $Y$, then the elementary theory of $\Gamma(Y)$ is decidable if and only if the elementary theory of $\Gamma(Y_i)$ is decidable for all $i=1,\dots, k$.
\item If $Y=Y_1\cup\dots \cup Y_k$ and $Z=Z_1\cup \dots \cup Z_l$ are two irreducible algebraic sets,  where $Y_1,\dots, Y_k$  and $Z_1,\dots, Z_l$ are the irreducible components of $Y$ and $Z$, respectively, then $\Gamma(Y)$ is elementary equivalent to $\Gamma(Z)$ if and only if $k=l$ and, after a certain re-enumeration, $\Gamma(Y_i)$ is elementary equivalent to $\Gamma(Z_i)$ for all $i=1,\dots, k$.
\end{enumerate}
\end{cor}

\section{Normal Forms of First-Order Formulas} \label{sec:mr}
\subsection{Conjunctions of positive formulas}

Let $L_{a,b}$ be the language of groups enriched by two constants $a$ and $b$, and let ${\mathcal S}$ be the class of all groups $G$ satisfying the following universal sentences:

\begin{itemize}
    \item $(I)\ \forall x \left( ([x,a]=1\wedge [x,b]=1) \to x=1 \right)$;
    \item $(II)\ \forall x \forall y \forall z (x^2 y^2 z^2 = 1 \rightarrow [x,y] = 1\ \wedge\ [x,z] = 1 \ \wedge \ [y,z] = 1)$;
    \item $(III)\ \forall x \forall y \left( x^2=y^2 \to x=y\right)$;
    \item $(IV)\ \forall x \left( [x^2, a]=1 \to [x,a]=1 \right)$.
\end{itemize}

Let GROUPS be a set of axioms of group theory. Denote by $A_{\mathcal S}$ the union of axioms $(I),(II),(III), (IV)$ and
GROUPS. Notice that the axiom $(II)$ is equivalent modulo GROUPS to the following quasi-identity
$$
\forall x \forall y \forall z (x^2 y^2 z^2 = 1 \rightarrow [x,y]=1).
$$
It follows that all axioms in $A_{\mathcal S}$ are quasi-identities.

\begin{lemma} \label{lem:pcinS}
The class $\mathcal S$ contains all partially commutative groups with trivial centre.
\end{lemma}

\begin{proof}

We first prove that in any partially commutative group $\GG$ with trivial centre there exist two elements $a$ and $b$ such that
$C(a)\cap C(b)=1$. Indeed, let $\GG=\GG_1\times \dots \times \GG_k$, be the decomposition of $\GG$ in the form (\ref{eq:decomp}). Since $Z(\GG)=1$, each $\GG_i$, $i=1,\dots, k$ is a non-abelian directly indecomposable partially commutative group. For each $i$ choose a pair of block elements $a_i, b_i\in \GG_i$ such that $C_{\GG_i}(a_i)\cap C_{\GG_i}(b_i)=1$. By Theorem \ref{thm:centr},
it follows that $C_\GG(a_1\dots a_k)=\langle \sqrt{a_1} \rangle \times \dots\times \langle \sqrt{a_k} \rangle$, $C_\GG(b_1\dots
b_k)=\langle \sqrt{b_1} \rangle \times \dots\times \langle \sqrt{b_k} \rangle$ and so $C_\GG(a_1\dots a_k)\cap C_\GG(b_1\dots
b_k)=1$. This proves that $\GG$ satisfies Axiom (I).

In \cite{crw} Crisp and Wiest prove the following theorem.

\begin{theorem}[J. Crisp, B. Wiest, \cite{crw}]
Let $G$ be a partially commutative group. Then the equation $x^2y^2z^2=1$ has only commutative solutions.
\end{theorem}

So $\GG$ satisfies Axiom (II).

By \cite{DK}, partially commutative groups have least roots, and thus $\GG$ satisfies Axiom (III).

By Corollary \ref{cor:centr}, $\GG$ satisfies Axiom (IV).
\end{proof}

\begin{lemma} \label{le:malcev2}
Let $G\in \mathcal S$. Then the equation
 \beq \label{eq:malcev}
    x^2ax^2a^{- 1} (ybyb^{-1})^{-2}=1
 \eeq
has only the trivial solution $x = 1$ and $y = 1$ in $G$.
\end{lemma}

\begin{proof}
Let $x,y$ be a solution of Equation (\ref{eq:malcev}) in $G$. Then we can rewrite (\ref{eq:malcev}) as follows
 \beq \label{eq:malcev2}
 (x^2 a)^2 a^{-2} = ((yb)^2 b^{-2})^2.
 \eeq
Since $G$ satisfies $(II)$, from (\ref{eq:malcev2}) we deduce that $[x^2a,a^{-1}]= 1$, hence $[x^2,a^{-1}]= 1$. Since $G$ satisfies $(IV)$,  it follows that $[x,a] = 1$. Now, we can rewrite (\ref{eq:malcev2}) in the form
$$
{(x^2)}^2 = ((yb)^2b^{-2})^2,
$$
and then, since $G$ satisfies $(III)$ we get
 \beq \label{eq:malcev3}
x^2 = (yb)^2b^{-2}.
 \eeq
Again, since $G$ satisfies $(II)$ it follows that $[x,b] = 1$ and $[y,b]=1$. This implies that $[x,a] = 1$ and $[x,b] =1$. Therefore,
applying $(I)$, we get $x=1$. In this event, (\ref{eq:malcev3}) reduces to $y^2 = 1$, so $y = 1$, as desired.
\end{proof}

\begin{cor}
\label{co:1} For any finite system of equations $S_1(X) = 1, \ldots, S_k(X) = 1$  one can effectively find a single equation $S(X) = 1$ such that given a group $G \in \mathcal{S}$, the following holds:
$$
V_G(S_1, \ldots, S_n) = V_G(S).
$$
\end{cor}
\begin{proof}
By induction it suffices to prove the result for $k = 2$. In this case, by the lemma above, the following equation
$$
S_1(X)^2aS_1(X)^2 a^{-1}(S_2(X)bS_2(X)b^{-1})^{-2}=1
$$
can be chosen as the equation $S(X) = 1.$
\end{proof}

\begin{cor}\label{co:2} For any finite system of atomic formulas
$$
S_1(X) = 1, \ldots, S_k(X) = 1
$$
in $L_{a,b}$, one can effectively find a atomic formula $S(X) = 1$ in $L_{a,b}$ such that $( \bigwedge\limits_{i = 1}^{k} S_i(X) = 1)$ is $A_{\mathcal S}$-equivalent to $S(X)=1$,
$$
\left( \bigwedge\limits_{i = 1}^{k} S_i(X) = 1\right)\ \sim_{A_{\mathcal S}} \ S(X) = 1.
$$
\end{cor}

\subsection{Disjunctions of positive formulas}
Our next aim is to be able to rewrite finite disjunctions of equations into conjunctions of equations.

Let ${\mathcal T}_D$ be the elementary theory (in the language $L_{a,b}$ of groups enriched by two constants) of non-abelian directly indecomposable partially commutative groups whose centraliser dimension is lower than a fixed number $D$, i.e.  the set of all first order sentences in the language of groups enriched by two constants $a$ and $b$ which are true in all non-abelian directly indecomposable partially commutative groups of centraliser dimension lower than $D$, together with the following two formulas:
\begin{itemize}
    \item The intersection of centralisers of $a$ and $b$ is trivial:
        $$
        \forall x ([x,a]=1 \wedge [x,b]=1) \rightarrow x=1
        $$
    \item The centralisers of $a$ and $b$ are cyclic. An interested reader may verify that this condition can indeed be written using the first order language.
\end{itemize}

\begin{rem}
 Note that any model of $\mathcal T _D$ lies in $\mathcal S$.
\end{rem}

In Section \ref{sec:posth}, we consider models in the language $L_\GG$ - the language of $\GG$-groups enriched by all constants from $\GG$ (where $\GG$ is a directly indecomposable, non-abelian partially commutative group). We denote the  elementary theory in the language $L_{\GG}$ of directly indecomposable partially commutative $\GG$-groups whose centraliser dimension is lower than a fixed number $D$ by the same symbol ${\mathcal T}_D$. The results of this and the next sections hold for both definitions of ${\mathcal T}_D$.

\begin{prop} \label{prop:66}
Let $G$ be a model of $\mathcal T_D$. Let $a, b\in G$ be elements with cyclic centralisers and such that $C(a)\cap C(b)=1$. Then for any solution $x,y \in G$ of the system
$$
[x,y] = 1, \ [x,y^{(a^{3\cdim(G)+4})}] = 1, \ \ [x,y^{(b^{3\cdim(G)+4})}] = 1,
$$
either $x = 1$ or $y = 1$.
\end{prop}
\begin{proof}
By Corollary \ref{cor:disj}, any directly indecomposable partially commutative group satisfies the following sentence in $L_{a,b}$:
$$
\forall x \forall y  ([x,y] = 1 \wedge [x,y^{(a^{3\cdim(\GG)+4})}] = 1 \wedge [x,y^{(b^{3\cdim(\GG)+4})}] = 1)\to (x = 1 \vee y = 1).
$$
Since the class of all groups that have centraliser dimension $D$ is universally axiomatisable (see \cite{DKR1}) any model $G$ of the theory $\mathcal T_D$ satisfies the above sentence and the statement follows.
\end{proof}

Combining Proposition \ref{prop:66} and Lemmas \ref{lem:pcinS} and \ref{le:malcev2} yields an algorithm to encode an arbitrary finite disjunction of equations into a single equation.

\begin{cor}
For any finite set of equations $S_1(X) = 1, \ldots, S_k(X) = 1$ one can effectively find a single equation $S(X) = 1$ such that given any model $G$ of ${\mathcal T}_D$, the following holds:
$$
V_G(S_1) \cup \ldots \cup V_G(S_k) = V_G(S).
$$
\end{cor}

\begin{cor}\label{co:4}
For any finite set of atomic formulas $S_1(X) = 1, \ldots, S_k(X) = 1$ one can effectively find a single atomic formula $S(X) = 1$ such that
$$
(\bigvee\limits_{i = 1}^{k} S_i(X) = 1) \ \sim_{{\mathcal T}_D} \ S(X) = 1.
$$
\end{cor}

\begin{cor} \label{co:qfreepos}
Every positive quantifier-free formula $\Phi(X)$ is equivalent modulo ${\mathcal T}_D$ to a single equation $S(X) = 1$.
\end{cor}

\subsection{Conjunctions and Disjunctions of Inequations}

The next result shows that one can effectively encode finite conjunctions and finite disjunctions of {\em inequations} (negations of atomic formulas) into a single inequation modulo ${\mathcal T}_D$.

\begin{lemma}
\label{le:ineq} For any finite set of inequations
$$
S_1(X) \neq 1, \ldots, S_k(X) \neq 1,
$$
one can effectively find an inequation $R(X) \neq 1$ and an inequation $T(X) \neq 1$ such that
$$
(\bigwedge_{i = 1}^{k} S_i(X) \neq 1)\ \sim_{{\mathcal T}_D} \ R(X) \neq 1
$$
and
$$
 (\bigvee_{i = 1}^{k}S_i(X) \neq 1) \ \sim_{{\mathcal T}_D} \ T(X) \neq 1.
$$
\end{lemma}
\begin{proof}
By Corollary \ref{co:4} there exists an equation $R(X) =  1$ such that
$$
\bigvee\limits_{i = 1}^{k} (S_i(X) = 1) \ \sim_{{\mathcal T}_D} \  R(X) = 1.
$$
Hence
$$
\left(\bigwedge\limits_{i =1}^{k} S_i(X) \neq 1 \right)  \  \sim_{{\mathcal T}_D}   \ \neg\left(\bigvee
\limits_{i = 1}^{k} S_i(X) = 1\right)  \ \sim_{{\mathcal T}_D} \ \neg (R(X) = 1)
\ \sim_{{\mathcal T}_D} \ R(X) \neq 1.
$$
This proves the first part of the result. Similarly, by Corollary \ref{co:2} there exists an equation $T(X) = 1$ such that
$$
\left( \bigwedge\limits_{i = 1}^{k} S_i(X) = 1\right) \  \sim_{{\mathcal T}_D} \ T(X)= 1.
$$
Hence
$$
\left( \bigvee\limits_{i = 1}^{k}S_i(X) \neq 1\right) \   \sim_{{\mathcal T}_D} \ \neg \left(\bigwedge\limits_{i = 1}^{k}S_i(X) = 1\right)
\  \sim_{{\mathcal T}_D} \ \neg (T(X) = 1) \  \sim_{{\mathcal T}_D}  \ T(X) \neq 1.
$$
\end{proof}

\begin{cor} \label{co:atomic}
For every quantifier-free formula $\Phi(X)$, one can effectively find a formula
$$
\Psi(X) = \bigvee\limits_{i=1}^n (S_i(X) = 1 \ \wedge \  T_i(X) \neq 1)
$$
which is equivalent to $\Phi(X)$ modulo ${\mathcal T}_D$. In particular, if $G$ is a model of ${\mathcal T}_D$, then  every quantifier-free formula $\Phi(X)$ is equivalent over $G$ to a formula  $\Psi(X)$ as above.
\end{cor}

\section{Positive Theory of Partially Commutative Groups}
\label{sec:posth}
In this section we present a procedure of quantifier elimination for positive formulas over partially commutative groups (an analog of Merzlyakov's Theorem for free groups). Our approach to the positive theory of partially commutative groups is based on the proof of Merzlyakov's Theorem given in \cite{IFT}.

\subsection{Generalised equations} \label{se:4-1}

Let $A= \{a_1, \ldots, a_m\}$ be a set of constants and $X = \{x_1, \ldots, x_n\}$ be a set of variables. Set  $\GG = G(A)$ to be a partially commutative group generated by $A$ and  $\GG[X] = \GG \ast F(X).$

\begin{defn}
A combinatorial generalised equation $\Omega$ (with constants from $A^{\pm 1}$) consists of the following objects:
\begin{enumerate}
    \item A finite set of {\em bases} $BS = BS(\Omega)$.  Every base is either a constant base or a variable base. Each constant base is associated with exactly one letter from $A^{\pm 1}$. The set of variable bases ${\mathcal M}$ consists of $2n$ elements ${\mathcal M} = \{\mu_1, \ldots, \mu_{2n}\}$. The set ${\mathcal M}$ comes equipped with two functions: a function $\varepsilon: {\mathcal M} \rightarrow \{1,-1\}$ and an involution $\Delta: {\mathcal M} \rightarrow {\mathcal M}$ (i.e., $\Delta$ is a bijection such that $\Delta^2$ is an identity on  ${\mathcal M}$). Bases $\mu$ and $\Delta(\mu)$ (or $\bar\mu$) are called {\em dual bases}.  We denote variable bases by $\mu, \lambda, \ldots.$
    \item A set of {\em boundaries} $BD = BD(\Omega)$. $BD$ is  a finite initial segment of the set of positive integers  $BD = \{1, 2, \ldots, \rho+1\}$. We use letters $i,j, \ldots$ for boundaries.
    \item Two functions $\alpha : BS \rightarrow BD$ and $\beta : BS \rightarrow BD$. We call $\alpha(\mu)$ and $\beta(\mu)$ the initial and terminal boundaries of the base $\mu$ (or endpoints of $\mu$). These functions satisfy the following conditions: $\alpha(b) <  \beta(b)$  for every base $b \in BS$; if $b$ is a constant base then $\beta(b) = \alpha(b) + 1$.
    \item A subset $C$ of $BD(\Omega)\times BD(\Omega)$.
\end{enumerate}
\end{defn}

To a combinatorial generalised equation $\Omega$ one can associate a system of equations in {\em variables} $h_1, \ldots, h_\rho$ over $G(A)$ (variables $h_i$ are sometimes  called {\em items}). This system is called a {\em generalised equation}, and, abusing the notation, we  denote it by the same symbol $\Omega$. The generalised equation  $\Omega$  consists of the
following three types of equations.

\begin{enumerate}
    \item Each pair of dual variable bases $(\lambda, \Delta(\lambda))$ provides an equation over a partially commutative group $\GG$
$$
[h_{\alpha (\lambda )}h_{\alpha (\lambda )+1}\cdots h_{\beta (\lambda )-1}]^ {\varepsilon (\lambda)}= [h_{\alpha (\Delta (\lambda ))}h_{\alpha (\Delta (\lambda ))+1} \cdots h_{\beta (\Delta (\lambda ))-1}]^ {\varepsilon (\Delta (\lambda))}.
$$
These equations are called {\em basic equations}. In the case when $\beta(\lambda)=\alpha(\lambda)+1$ and $\beta(\Delta(\lambda))=\alpha(\Delta(\lambda))+1$, i.e. the corresponding basic equation takes the form:
$$
[h_{\alpha (\lambda )}]^ {\varepsilon (\lambda)}= [h_{\alpha (\Delta (\lambda ))}]^ {\varepsilon (\Delta (\lambda))},
$$
without loss of generality, we shall assume that the equality above is graphical.

    \item For each constant base $b$ we write down a {\em coefficient equation}
$$
h_{\alpha(b)} = a,
$$
where $a \in A^{\pm 1}$ is the constant associated with $b$.

    \item For every element $c=(i,j)\in C$ we write the following equation:
    $$
    [h_i,h_j]=1.
    $$
\end{enumerate}

\begin{rem}
We assume that every generalised equation comes associated with a combinatorial one;
\end{rem}

Let $\GG=G(A)$, then the monoid given by the presentation
$$
\Tr=\Tr(A^{\pm 1})= \left< A\cup A^{-1}\mid [ a_i^{\pm 1}, a_j^{\pm 1} ] =1 \right>, \hbox{ where } \left\{ \left[ a_i^{\pm 1}, a_j^{\pm 1} \right] =1 \hbox{ in } \GG \right\}
$$
is called \emph{partially commutative monoid associated to} $\GG$. Partially commutative monoids, are also known as \emph{trace} monoids and are extensively studied, see \cite{traces} and references there.

\begin{defn}
Let $\Omega(h) = \{L_1(h)=R_1(h), \ldots, L_s(h) = R_s(h)\}$ be a generalised equation in variables $h = (h_1, \ldots,h_{\rho})$ with constants from $A^{\pm 1 }$. A sequence of reduced nonempty words $U = (U_1(A), \ldots, U_{\rho}(A))$ in the alphabet $A^{\pm 1}$ is a {\em solution} of $\Omega $ if:
\begin{enumerate}
\item [1)] all words $L_i(U), R_i(U)$ are geodesic (treated as elements of $\GG$) as written;

\item [2)] $L_i(U) =  R_i(U),  \quad i = 1, \ldots s$ in the partially commutative monoid $\Tr(A^{\pm 1})$.
\end{enumerate}
\end{defn}

The notation $(\Omega, U)$ means that $U$ is a solution of the generalised equation $\Omega $.

\begin{rem}
Notice that a solution $U$ of a generalised equation $\Omega$ can be viewed as  a solution of $\Omega$ in the partially commutative monoid $\Tr(A^{\pm 1})$ (i.e., $L_i(U) = R_i(U)$ modulo commutation) which satisfies an additional condition:  $U \in {\Tr(A^{\pm 1})}^\rho$ and $U$ is a tuple of geodesic words in $\GG$.
\end{rem}

Obviously, each solution  $U$ of $\Omega$ gives rise to  a solution of $\Omega$ in the partially commutative group $G(A)$. The converse does not hold in general, i.e. it might happen that $U$ is a solution of $\Omega$ in $G(A)$ but not in $\Tr(A^{\pm 1})$, i.e. some equalities $L_i(U) = R_i(U)$ hold only after a reduction in $\GG$. We introduce the following notation which will allow us  to distinguish in which structure ($\Tr(A^{\pm 1})$ or $G(A)$) we are resolving $\Omega$.

If
$$
S = \{L_1(h)=R_1(h), \ldots, L_s(h) = R_s(h)\}
$$
is an arbitrary system of equations with constants from $A^{\pm 1}$, then by $S^*$ we denote the system of equations
$$
S^*= \{L_1(h)R_1(h)^{-1} = 1, \ldots, L_s(h)R_s(h)^{-1} = 1\}
$$
over the group $G(A)$.

\subsection{Reduction to generalised equations} \label{sec:redge}

Similarly to the case of free groups, we now show how for a given finite system of equations $S(X,A) = 1$ over a partially commutative group $\GG$ one can associate a finite collection  of generalised equations $\mathcal{GE}(S)$ with constants from $A^{\pm 1}$. The collection $\mathcal{GE}(S)$ to some extent describes all solutions of the system $S(X,A) = 1$.

Informally, Lemma \ref{lem:prod} describes all possible cancellation schemes for the set of all solutions of the system $S(X,A)$ in the following way: the cancellation scheme corresponding to a particular solution, can be obtained from the one described in Lemma \ref{lem:prod} by setting some of the words $w_i^j$'s (and the corresponding bands) to be trivial. Therefore, every partition table (to be defined below) corresponds to one of the cancellation schemes obtained from the general one by setting some of the words $w_i^j$'s to be trivial. Every non-trivial word $w_i^j$ corresponds to a variable $z_k$ and the word $w_j^i$ to the variable $z_k^{-1}$.  If a variable $x$ that occurs in the system $S(X,A) = 1$ is subdivided as a product of some words $w_i^j$'s, i.e. is a word in the $w_i^j$'s, then the word $V_{ij}$ from the definition of a partition table is this word in the corresponding $z_k$'s. If the bands corresponding to the words $w_i^j$ and $w_k^l$ cross, then the corresponding variables $z_r$ and $z_s$ commute in the group $\Gamma$. We refer the reader to the construction of a partition table by a solution of the system $S(X,A)$ given in the end of this section to gain an intuition of the definition of a partition table.

Write $\{S(X,A)=1\}= \{S_ 1 = 1, \ldots, S_m = 1\}$ in the form
\begin{equation}\label{*}
\begin{array}{c}
 r_{11}r_{12}\ldots r_{1l_1}=1,\\
 r_{21}r_{22}\ldots r_{2l_2}=1,\\
 \ldots \\
 r_{m1}r_{m2}\ldots r_{ml_m}=1,\\
\end{array}
\end{equation}
where $r_{ij}$ are letters of the alphabet $ X^{\pm 1}\cup A^{\pm 1}$.

A pair (a set of geodesic words, a $\GG$-partially commutative group), $T=(V,\Gamma)$ of the form:
$$
V = \{V_{ij}(z_1, \ldots ,z_p)\}\subset \GG*F(Z)=\GG[Z] \ \ (1\leq i\leq m, 1\leq j\leq l_i), \quad \Gamma=\GG(A\cup Z),
$$
is called a {\em partition table} of the system $S(X,A)$ if the following conditions are satisfied:
\begin{enumerate}
    \item [0)] Every letter from $Z\cup Z^{-1}$ occurs in the words $V_{ij}$, moreover it occurs only once;
    \item [1)] The equality $V_{i1}V_{i2} \ldots V_{il_i}=1, 1\leq i\leq m,$ holds in $\Gamma$;
    \item [2)] $|V_{ij}|\leq l_i-1$;
    \item [3)] if $r_{ij} = a \in A^{\pm 1}$, then $|V_{ij}|=1$.
\end{enumerate}

Since $|V_{ij}|\leq l_i - 1$ then at most $|S(X,A)| = \sum\limits_{i = 1}^m (l_i - 1)l_i$ different letters $z_i$ can occur in a partition
table of $S(X,A) = 1$. Therefore we will always assume that $p \leq |S|$. We call the number $|S(X,A)|$ the \emph{size} of the system $S$.

Each partition table encodes a particular type of cancelation that happens when one substitutes a particular solution $W(A) \in G(A)$ into $S(X,A) = 1$ and then reduces (in a certain way) the words in $S(W(A),A)$ into the empty word.

\begin{lemma}\label{le:4.2}
Let $S(X,A) = 1$ be a finite system of equations over $G(A)$. Then
\begin{enumerate}
    \item the set $PT(S)$ of all partition tables  of $S(X,A) = 1$ is finite, and its cardinality is bounded by a number which depends only on  $|S(X,A)|$;
    \item one can effectively enumerate the set $PT(S)$.
\end{enumerate}
\end{lemma}
\begin{proof}   Since the words $V_{ij}$ have  bounded length, one can effectively enumerate the finite set  of all collections of words $\{V_{ij}\}$ in $\GG[Z]$ which satisfy the conditions 0), 2), 3) above.  Now for each such collection $\{V_{ij}\}$,  one can
effectively check  whether the  equalities $V_{i1}V_{i2} \ldots V_{il_i}=1, 1\leq i\leq m$ hold in one of the finitely many (since $|Z|<\infty$) partially commutative groups $\Gamma$ or not. This allows one to list effectively all partition tables for $S(X,A) = 1$.
\end{proof}

To each partition table $T=(\{V_{ij}\}, \Gamma)$ one can assign a generalised equation $\Omega _T$ in the following way (below we
use `$\doteq$' for graphical equality, i.e. equality in the free monoid). Consider the following word $V$ in $M(A^{\pm 1} \cup Z^{\pm 1})$:
$$
V\doteq V_{11}V_{12}\ldots V_{1l_1}\ldots V_{m1}V_{m2} \ldots V_{ml_m} = y_1 \ldots y_\rho,
$$
where $y_i \in A^{\pm 1} \cup Z^{\pm 1}$ and $\rho = l(V)$ is the length of $V$. Then the generalised equation $\Omega_T =
\Omega_T(h)$ has $\rho + 1$ boundaries and $\rho$ variables $h_1,\ldots ,h_{\rho}$ which are denoted by   $h = (h_1,\ldots, h_{\rho})$.

Now we define bases of $\Omega_T$ and the functions $\alpha, \beta, \varepsilon$.

Let $z \in Z$. For the (only) pair of occurrences of $z$ in $V$:
$$
y_i = z^{\epsilon _i}, \quad y_j = z^{\epsilon _j} \quad (\epsilon _i, \epsilon _j \in \{1,-1\})
$$
we introduce a pair of dual variable bases $\mu_{z,i}, \mu_{z,j}$ such that $\Delta(\mu_{z,i}) = \mu_{z,j}$. Put
\begin{gather}\notag
\begin{split}
\alpha(\mu_{z,i}) = i, \quad &\beta(\mu_{z,i}) = i+1, \quad \varepsilon(\mu_{z,i}) = \epsilon _i,\\
&\alpha(\mu_{z,j}) = j, \quad \beta(\mu_{z,j}) = j+1, \quad \varepsilon(\mu_{z,j}) = \epsilon _j.
\end{split}
\end{gather}
The basic equation that corresponds to this pair of dual bases is $h_{i}^{\epsilon_i}\doteq h_{j}^{\epsilon _j}$.

Let $x \in X$. For any two distinct occurrences of $x$ in $S(X,A) = 1$:
$$
 r_{i,j} = x^{\epsilon_{ij}}, \quad r_{s,t} = x^{\epsilon_{st}} \quad (\epsilon _{ij}, \epsilon _{st} \in \{1,-1\})
$$
so that $(i,j)$ precedes $(s,t)$ in left-lexicographical order, we introduce a pair of dual bases $\mu_{x,(ij),(st)}$ and $\mu_{x,(st),(ij)}$ such that $\Delta(\mu_{x,i,j,s,t}) = \mu_{x,s,t,i,j}$. Let $V_{ij}$ occur in the word $V$ as a subword
$$
V_{ij} = y_{c_1} \ldots y_{d_1}, \quad V_{st} = y_{c_2} \ldots y_{d_2}.
$$
Then we put
\begin{gather}\notag
\begin{split}
\alpha(\mu_{x,i,j,s,t}) = {c_1}, \quad &\beta(\mu_{x,i,j,s,t}) = d_1+1, \quad \varepsilon(\mu_{x,i,j,s,t}) = \epsilon_{ij},\\
&\alpha(\mu_{x,s,t,i,j}) = {c_2}, \quad \beta(\mu_{x,s,t,i,j}) = d_2+1, \quad \varepsilon(\mu_{x,s,t,i,j}) = \epsilon_{st}.
\end{split}
\end{gather}
The basic equation which corresponds to these dual bases can be written in the form
$$
[h_{\alpha(\mu_{x,i,j,s,t})} \ldots h_{\beta(\mu_{x,i,j,s,t})-1}]^{\epsilon _{ij}}=_{\Tr} [h_{\alpha(\mu_{x,s,t,i,j})}\ldots h_{\beta(\mu_{x,s,t,i,j})-1}]^{\epsilon _{st}}.
$$

Let $r_{ij} = a \in A^{\pm 1}$. In this case we introduce a constant base $\mu_{ij}$ with the label $a$. If $V_{ij}$ occurs in $V$ as $V_{ij} = y_c$, then we put
$$
\alpha(\mu_{ij}) = c, \quad \beta(\mu_{ij}) = c+1.
$$
The corresponding coefficient equation is written as $h_{c}=a$.

For any two distinct occurrences  $z_1, z_2 \in V$ such that  $[z_1,z_2]=1$ in $\Gamma$:
$$
y_i={z_1}^\epsilon_1, \quad y_j={z_2}^{\epsilon_j}, \quad (\epsilon_i, \epsilon_j\in \{1,-1\})
$$
we set $(i,j)\in C$. The corresponding equation is $[h_i,h_j]=1$. This defines the generalised equation $\Omega_T$.
Put
$$
\mathcal{GE}(S) = \{\Omega_T \mid T \hbox{ is a partition table for } S(X,A)= 1 \}.
$$
Then $\mathcal{GE}(S)$ is a finite collection of generalised equations which can be effectively constructed for a given $S(X,A) = 1$.

By $\GG_{R(\Omega)}$ we denote the coordinate group of $\Omega^*$, $\GG_{R(\Omega)}=\factor{\GG[h]}{R(\Omega^*)}$ (recall that by $\Omega^*$ we denote the system of equation over the group $G(A)$). Now we explain relations between  the coordinate groups of $S(X,A) = 1$ and of $\Omega_T^*$.

For a letter $x$ in $X$ we choose an arbitrary occurrence of $x$ in $S(X,A) = 1$ as
$$
r_{ij} = x^{\epsilon_{ij}}.
$$
Let $\mu = \mu_{x,i,j,s,t}$ be a base that corresponds to this occurrence of $x$. Then $V_{ij}$ occurs in $V$ as the subword
$$
V_{ij} = y_{\alpha(\mu)} \ldots y_{\beta(\mu) -1}.
$$
Notice that the word $V_{ij}$ does not depend on the choice of the base $\mu_{x,i,j,s,t}$  corresponding to the occurrence $r_{ij}$.

Define a word $P_x(h) \in \GG[h]$ (where $h = \{h_1, \ldots,h_\rho\}$) as follows
$$
P_x(h,A) =\left( h_{\alpha(\mu)} \ldots h_{\beta(\mu)-1}\right)^{\epsilon_{ij}},
$$
and put
$$
P(h) = (P_{x_1}, \ldots, P_{x_n}).
$$

The tuple of words $P(h)$ depends only on the choice of occurrences of letters from  $X$ in $V$.
It follows from the construction above that the map $X \rightarrow \GG[h]$ defined by $x \rightarrow  P_x(h,A)$ gives rise to a $\GG$-homomorphism
$$
\pi : \GG_{R(S)}\rightarrow \GG_{R(\Omega _T)}.
$$
Indeed,  if $f(X) \in R(S)$ then $\pi(f(X))= f(P(h))$. Then given a solution of $\Omega_T$ it follows from condition 1) of the definition of partition table that $f(P(h)) = 1$, thus $R(f(S))\subseteq R(\Omega^*_T)$ and $\pi$ is a homomorphism.

Observe that the image  $\pi (x)$ in $\GG_{R(\Omega _T)}$ does not depend on a particular choice of the occurrence of $x$ in $S(X,A)$ (the basic equations of $\Omega_T$ make these images equal). Hence $\pi$ depends only on $\Omega_T$.

To relate solutions of $S(X,A) = 1$  to solutions of generalised  equations  from $\mathcal{GE}(S)$ we need the technique developed in Section \ref{sec:cancan}.

Let $W(A)$ be a solution of $S(X,A) = 1$ in $G(A)$. If in the system (\ref{*}) we make the substitution  $\sigma : X \rightarrow W(A)$, then
$$
(r_{i1}r_{i2}\ldots r_{il_i})^{\sigma} = r_{i1}^{\sigma}r_{i2}^{\sigma}\ldots r_{il_i}^{\sigma} = 1
$$
in $G(A)$ for every $i = 1, \ldots, m$.

Since every product $R_i = r_{i1}^{\sigma}r_{i2}^{\sigma}\ldots r_{il_i}^{\sigma}$ is trivial, we can choose a van Kampen diagram $\D_{R_i}$ for $R_i$. Denote by ${\tilde z}_{i,1}, \ldots, {\tilde z}_{i,p_i}$ the subwords $w_j^k$, $1\le j<k\le l_i$ of $r_{ij}^{\sigma}$, where here $w_j^k$ are defined as in Lemma \ref{lem:prod}. Since, by Lemma \ref{lem:prod} $w_j^k= {w_k^j}^{-1}$, the word $r_{ij}$ can be written as a word in ${\tilde z}_{i,1}, \ldots,  {\tilde z}_{i,p_i}$:
$$
r_{ij}^\sigma = V_{ij}({\tilde z}_{i,1}, \ldots, {\tilde z}_{i,p_i})
$$
for some freely reduced words $V_{ij}(Z_i)$ in variables $Z_i = \{z_{i,1}, \ldots, z_{i,p_i}\}$. Observe that if $r_{ij} = a \in A^{\pm 1}$ then $r_{ij}^{\sigma} = a$ and we have $|V_{ij}| = 1$. By Lemma \ref{lem:prod}, $r_{ij}^\sigma$ is a product of at most $l_i-1$ words $w_j^k$, we have $|V_{ij}| \leq l_i - 1$. Denote by $Z=\bigcup\limits_{i=1}^m Z_i= \{z_1, \ldots, z_p\}$. Take a partially commutative group $\Gamma=G(A\cup Z)$ whose underlying commutation graph is defined as follows:
\begin{itemize}
    \item two elements $a_i,a_j$ in $A^{\pm 1}$ commute whenever they commute in $\GG$;
    \item an element $a\in A^{\pm 1}$ commutes with $z_i$ if and only if $a$ commutes with the word $w_j^k$ corresponding to $z_i$;
    \item two elements $z_i, z_j\in Z$ commute whenever the corresponding words $w_j^k$ do.
\end{itemize}

In the above notation, the set $T = (\{V_{ij}\}, \Gamma)$ is a partition table for $S(X,A) = 1$. We define
$$
U(A) = ({\tilde z}_1, \ldots, {\tilde z}_p)
$$
to be the solution of the generalised equation $\Omega_T$ induced by $W(A)$. From the construction of the map $P(h)$ we
deduce that $W(A) = P(U(A))$.

The converse is also true: if $U(A)$ is an arbitrary  solution of the generalised equation $\Omega_T$, then $P(U(A))$ is a solution of $S(X,A) = 1$.

We summarize the discussion above in the following lemma.

\begin{lemma}\label{le:R1}
For a given system of equations $S(X,A)=1$ over $\GG$, one can effectively construct a finite set
$$
\mathcal{GE}(S) = \{\Omega_T \mid T \hbox{ is a partition table for } S(X,A)= 1 \}
$$
of generalised equations such that
\begin{enumerate}
    \item if the set $\mathcal{GE}(S)$ is empty, then $S(X,A)= 1$ has no solutions in $\GG$;
    \item for each $\Omega (h) \in\mathcal{GE}(S)$ and for each $x \in X$ one can effectively find a word $P_x(h,A) \in \GG[h]$ of length at most $|h|$ such that the map $x\mapsto P_x(h,A)$ {\rm(}$x \in X${\rm)} gives rise to a $\GG$-homomorphism $\pi_\Omega: \GG_{R(S)}\rightarrow \GG_{R(\Omega)}$;
    \item for any solution $W(A) \in \GG^n$ of the system $S(X,A)=1$ there exists $\Omega (h) \in \mathcal{GE}(S)$ and a solution $U(A)$ of $\Omega(h)$ such that $W(A) = P(U(A))$, where $P(h) = (P_{x_1}, \ldots, P_{x_n})$, and this equality holds in the partially commutative monoid $\Tr(A^{\pm 1})$;
    \item for any $\GG$-group $\tilde G$, if a  generalised equation $\Omega (h) \in \mathcal{GE}(S)$  has a solution $\tilde U$ in $\tilde G$, then $P(\tilde U)$ is a solution of $S(X,A) = 1$ in $\tilde G$.
\end{enumerate}
\end{lemma}

\begin{cor} \label{co:R1}
In the notation of {\rm Lemma \ref{le:R1}}  for any solution $W(A)  \in \GG^n=G(A)^n$ of the system $S(X,A)=1$ there exists $\Omega (h) \in \mathcal{GE}(S)$  and a solution $U(A)$ of $\Omega(h)$ such that the following diagram commutes
$$
\xymatrix@C3em{
 \GG_{R(S)}  \ar[rd]_{\pi_W} \ar[rr]^{\pi}  &  &\GG_{R(\Omega)} \ar[ld]^{\pi_U}
                                                                             \\
                               &  \GG &
}
$$
\end{cor}

\subsection{Positive theory of partially commutative groups and direct products of groups}\label{se:mer}
In this section we first prove a result on elimination of quantifiers for positive sentences over a non-abelian directly indecomposable partially commutative group $\GG = G(A)$. This proof is based on the notion of a generalised equation. Combining this result with a theorem of V.~Diekert and A.~Muscholl on decidability of equations over partially commutative groups, see \cite{DM}, we get that the positive theory of free partially commutative groups in the language of group theory $L$ and the language $L_\GG$ enriched by constants is decidable. Note that V.~Diekert and M.~Lohrey, using a different method, prove a similar result in \cite{DL}. Furthermore, we apply the techniques developed for the proof of quantifier elimination to obtain a result on lifting arbitrary formulas from $\GG$ to $\GG\ast F$, where $F$ is a free group of finite rank (see Theorem \ref{merzl}). In order to prove that the positive theory of any partially commutative group is decidable, we need to study the positive theory of the direct product of groups. In the appendix of the paper we prove that if $G=H_1\times \cdots \times H_k$, then the positive theory of $G$ in both languages $L$ and $L_G$ is decidable if the positive theories of $H_1, \dots, H_k$ are decidable in $L$ and $L_{H_i}$, correspondingly.

Recall that every positive formula $\Psi(Z)$ in the language $L_G$ is equivalent modulo ${\mathcal T}_D$ to a formula of the type
$$
\forall x_1 \exists y_1 \ldots \forall x_k \exists y_k (S(X,Y,Z,A)= 1),
$$
where $S(X,Y,Z,A) = 1$ is an equation with constants from $A^{\pm 1}$, $X = (x_1, \ldots, x_k)$, $Y = (y_1, \ldots,y_k)$, $Z =
(z_1, \ldots, z_m)$. Indeed, one can insert auxiliary quantifiers to ensure the direct alteration of quantifiers in the prefix.
In particular, every positive sentence in $L_G$ is equivalent modulo ${\mathcal T}_D$ to a formula of the type
$$
\forall x_1 \exists y_1 \ldots \forall x_k \exists y_k (S(X,Y,A)= 1).
$$

\begin{theorem}[Elimination of Quantifiers] \label{thm:eq}
If
$$
\GG \models \forall x_1\exists y_1\ldots\forall x_k\exists y_k (S(X,Y,A)=1),
$$
then there exist words {\rm (}with constants from $\GG${\rm)} $q_1(x_1),\ldots , q_k(x_1,\ldots ,x_k) \in \GG[X],$ such that
$$
\GG[X]\models S(x_1, q_1(x_1),\ldots ,x_k, q_k(x_1,\ldots,x_k),A)=1,
$$
i.e. the equation
$$
S(x_1,y_1, \ldots, x_k,y_k,A) = 1
$$
{\rm(}in variables $Y${\rm)} has a solution in the group $\GG[X]$.
\end{theorem}
\begin{proof}
Let $\mathcal{GE}(S) = \{\Omega_1(Z_1), \ldots, \Omega_r(Z_r)\}$ be generalised equations associated with the equation $S(X,Y,A) = 1$ in Lemma \ref{le:R1}. Denote by  $\rho_i = |Z_i|$ the number of variables in $\Omega_i$.

Since the group $\GG$ is directly indecomposable, there exists a path $p$ in the non-commutation graph $\Delta$ of $\GG$ beginning in a vertex $b_1$ which goes through every vertex of $\Delta$ at least once. Denote by $b_1\cdots b_n$ the label of the path $p$. Set
$$
b=b_1b_2\cdots b_{n-1}b_nb_{n-1}\cdots b_1, \quad a=b_2bb_2=b_2b_1b_2\cdots b_{n-1}b_nb_{n-1}\cdots b_2 b_1 b_2,
$$
and
$$
g_{1}=b^ma^{m_{1,1}}b^m a^{m_{1,2}}b\ldots a^{m_{1,n_{1}}}b^m,
$$
where $0<m_{1,1} < m_{1,2} < \ldots < m_{1,n_1}$, $\max\{\rho_1, \ldots,\rho_r\}|S(X,A)| < n_1$ (recall that by $|S(X,A)|$ we denote the size of the system $S$, see Section \ref{sec:redge})  and $m \in \BN$ is a constant which depends only on the generalised equation.

For the word $g_1$ there exists $h_1$ such that
$$
\GG \models \forall x_2 \exists y_2\ldots\forall x_k\exists y_k (S(g_1,h_1, x_2, y_2, \ldots, x_k,y_k)=1).
$$
Suppose now that elements $g_1, h_1, \ldots g_{i- 1},h_{i-1} \in \GG$ are so that
$$
\GG \models \forall x_i \exists y_i\ldots\forall x_k\exists y_k (S(g_1,h_1, \ldots, g_{i-1},h_{i-1},x_i, y_i, \ldots, x_k,y_k)=1).
$$
We define
\begin{equation}
\label{eq:two-stars}
g_{i} = b^m a^{m_{i1}}b^m a^{m_{i2}}b^m\ldots a^{m_{in_{i}}}b^m
\end{equation}
such that
\begin{enumerate}
    \item [1)] $0<m_{i1} < m_{i2} < \ldots m_{in_i}$;
    \item [2)] $\max\{\rho_1, \ldots,\rho_r\} |S(X,A)| < n_i $;
    \item [3)] no subword of the type $b^m a^{m_{ij}} b^m$ occurs in any of the words $g_l$, $l < i$ and in any of the (finitely many) words $h'_l$ such that $h'_l=h_l$ in  $\Tr(A^{\pm 1})$, $l < i$.
\end{enumerate}
Then there exists an element $h_i \in \GG$ such that
$$
\GG \models \forall x_{i+1}\exists y_{i+1}\ldots\forall x_k\exists y_k  (S(g_1,h_1, \ldots,  g_i, h_i,x_{i+1},y_{i+1},  \ldots, x_k,y_k) =1).
$$
By induction we have constructed elements $g_1, h_1, \ldots, g_k, h_k \in \GG$ such that
$$
S(g_1,h_1, \ldots,g_k,h_k) = 1
$$
and each $g_i$ has the form (\ref{eq:two-stars}) and satisfies conditions 1), 2), 3).

By Lemma \ref{le:R1}, there exists a generalised equation $\Omega(Z) \in \mathcal{GE}(S)$, words $P_{i}(Z,A), Q_{i}(Z,A)
\in \GG[Z]$ ($i = 1,\ldots,k$) of length lower than $\rho = |Z|$, and a solution $U = (u_1, \ldots, u_\rho)$ of $\Omega(Z)$ in
$\GG$ such that the following words are equal in $\Tr(A^{\pm 1})$:
$$
g_i = P_{i}(U), \quad h_i = Q_{i}(U) \quad (i = 1, \ldots,k).
$$
Notice that from the definition of $a$ and $b$ it follows that no two consecutive letters in $a$ and $b$, and thus in $g_i$ commute. Therefore, the equality $g_i = P_{i}(U)$ is graphical (i.e. $g_i = P_{i}(U)$ in the free monoid).

Since $n_i > \rho|S(X,A)|$ (by condition 2)) and $P_i(U) = y_1 \ldots y_q$ with  $y_i \in U^{\pm 1}, q \leq \rho$, the graphical  equalities
\begin{equation} \label{eq:gp}
g_i = b^m a^{m_{i1}}b^m a^{m_{i2}}b^m \ldots a^{m_{in_{i}}}b^m = P_{i}(U) \quad (i = 1,\ldots,k)
\end{equation}
show that  there exists a subword $v_i = b^m a^{m_{il}} b^m$ of $g_i$ such that every occurrence of this subword in (\ref{eq:gp})  is an occurrence inside some $u_j^{\pm 1}$. For each $i$ fix  such a subword  $v_i = b^m a^{m_{il}} b^m$ in $g_i$. In view of condition 3), the word $v_i$ does not occur in any of the words $g_j$ ($j \neq i$), $h_s$ ($s < i$), moreover, in $g_i$ it occurs precisely once. Denote by $j(i)$ the unique index such that $v_i$ occurs inside $u_{j(i)}^{\pm 1}$ in $P_i(U)$ from (\ref{eq:gp}) (and $v_i$ occurs in it precisely once).

The argument above shows that the variable $z_{j(i)}$ does not occur in words $P_{t}(Z,A)$ ($t \neq i$), $Q_{s}(Z,A)$ ($s < i$).
Moreover, in $P_{i}(Z)$ it occurs precisely once.  It follows that the variable $z_{j(i)}$ in the generalised equation $\Omega(Z)$ does not occur in  either the coefficient equations or in the basic equations corresponding to the dual bases related to $x_t$ ($t \neq i$), $y_s$ ($s < i$).

We ``mark'' (or select) the unique occurrence of $v_i$ (as $v_i^{\pm 1}$) in $u_{j(i)}$ $i = 1, \ldots, k$. Now we are going
to mark some other occurrences of $v_i$ in the words $u_1, \ldots, u_\rho$ as follows. Suppose that some $u_d$ has a marked occurrence of some $v_i$. If $\Omega$ contains an equation of the type $z_d^{\epsilon} =  z_r^{\delta}$, then $u_d^{\epsilon}= u_r^{\delta}$ graphically. Hence $u_r$ has an occurrence of the subword $v_i^{\pm 1}$ which corresponds to the marked occurrence of $v_i^{\pm 1}$ in $u_d$. We mark this occurrence of $v_i^{\pm 1}$ in $u_r$.

Suppose $\Omega$ contains an equation of the type
$$
[z_{\alpha _1}\ldots z_{\beta _1-1}]^{\epsilon _1}= [z_{\alpha _2}\ldots z_{\beta _2- 1}]^{\epsilon _2}
$$
such that $z_d$ occurs in it, say in the left-hand side of the equality. Then
$$
[u_{\alpha _1}\ldots u_{\beta_1-1}]^{\epsilon _1}= [u_{\alpha _2}\ldots u_{\beta _2-1}]^{\epsilon _2}
$$
in the partially commutative monoid $\Tr(A^{\pm 1})$. Since  $v_i^{\pm 1}$ is a subword of $u_d$, a subword $v_{i,1}=b^{m-1}a^{m_{il}}b^{m-1}$ occurs also in the right-hand side of the above equality, say in $u_r$. Indeed, let $w_1 b^{m}a^{m_{il}}b^{m}w_2=w$ in the monoid $\Tr(A^{\pm 1})$. Since for any letter $\ell$ in $w_k$, $k=1,2$ there exists a letter $\ell'$ in $b$ such that $[\ell,\ell'] \ne 1$ and since, by the definition of $a$ and $b$, no two consecutive occurrences in $b^{m-1}a^{m_{il}}b^{m-1}$ commute, the statement follows. We mark this occurrence of $v_{i,1}^{\pm 1}$ in  $u_r$ and in all the previously marked occurrences of $v_i^{\pm 1}=bv_{i,1}b$. We continue the marking process, but now, instead of $v_i$ we mark the occurrences of $v_{i,1}$. The marking process stops in finitely many steps and all the occurrences of the subword $v_{i,k}=b^{m-k}a^{m_{il}}b^{m-k}$ are marked. For the above argument, it suffices to choose $m>k$, which depends on the generalised equation only.

Now in all words $u_1, \ldots, u_\rho$ we replace every marked occurrence of $v_{i,k} = b^{m-k}a^{m_{il}}b^{m-k}$ with a new word $b^{m-k}a^{m_{il}}x_{i}b^{m-k}$ from the group $\GG[X]$. Denote the resulting words from $\GG[X]$ by $\tilde u_1, \ldots, \tilde u_\rho$. It follows from description of the marking process that the tuple $\tilde U = (\tilde u_1, \ldots, \tilde u_{\rho})$  is a solution of the generalised equation $\Omega$ in $\GG[X]$. Indeed, by construction, all basic and coefficient equations in $\Omega$ hold in the partially commutative monoid if we substitute $z_i$ by $\tilde u_i$. Furthermore, since any word that contains the word of the form $b^{m-k}a^{m_{il}}x_{i}b^{m-k}$ has cyclic centraliser, it follows that $\tilde U$ satisfies the same commutation equations as $U$ and thus the commutation equations of $\Omega$. Now, by Lemma \ref{le:R1}, $X = P(\tilde U), Y = Q(\tilde U)$ is a solution of the equation $S(X,A) = 1$ over  $\GG[X]$ as desired.
\end{proof}

\begin{cor} \label{co:merz}
There is an algorithm which for a given positive sentence
$$
\forall x_1\exists y_1\ldots\forall x_k\exists y_k (S(X,Y,A)=1)
$$
in $L_\GG$ determines whether or not this formula holds in $\GG$, and if it does, the algorithm finds words
$$
q_1(x_1),\ldots , q_k(x_1,\ldots ,x_k) \in \GG[X]
$$
such that
$$
\GG[X]\models S(x_1, q_1(x_1),\ldots ,x_k, q_k(x_1,\ldots ,x_k),A)=1,
$$
i.e. the positive theory of any directly-indecomposable partially commutative group is decidable.
\end{cor}
\begin{proof}   The proof follows from Theorem \ref{thm:eq} and decidability of equations over free partially commutative groups. Indeed, the compatibility problem for a system of equations over a partially commutative group $\GG$ reduces to the compatibility of a system of equations $S'$ over a free group with rational constraints $C$, see \cite{DM}. In order to prove the Corollary it suffices to check the compatibility of $S'$ over a free group with constraints $C \cup \{y_i \in F[X_i]\}$, where $X_i = \{x_1, \ldots, x_i\}$.
\end{proof}

The next result follows directly from Corollary \ref{cor:pthdir}.
\begin{cor}
Let $\GG$ be an arbitrary partially commutative group. Then the positive theory of $\GG$ is decidable.
\end{cor}

\begin{defn}
Let $\phi$ be a sentence in the language $L_\GG$ written in the standard form
$$
\phi = \forall x_1\exists y_1\ldots\forall x_k\exists y_k \ \phi_0(x_1,y_1, \ldots, x_k,y_k),
$$
where $\phi_0$ is a quantifier-free formula in $L_\GG$. We say that $\GG$ {\em freely lifts} $\phi$ if there exist words
(with constants from $\GG$) $q_1(x_1),\ldots , q_k(x_1,\ldots ,x_k) \in \GG[X],$ such that
$$
\GG[X]\models \ \phi_0(x_1, q_1(x_1),\ldots ,x_k, q_k(x_1,\ldots,x_k),A)=1.
$$
\end{defn}

\begin{theorem} \label{merzl}
Let $\GG$ be a non-abelian directly indecomposable partially commutative group. Then $\GG$ freely lifts every sentence in $L_\GG$ that is true in $\GG$.
\end{theorem}
\begin{proof}
Suppose a sentence
\begin{equation} \label{eq:sent}
\phi = \forall x_1\exists y_1\ldots\forall x_k\exists y_k (U(x_1,y_1, \ldots, x_k,y_k) =1  \wedge V(x_1,y_1, \ldots, x_k,y_k) \neq 1),
\end{equation}
is true in $\GG$. We choose $x_1 = g_1, y_1 = h_1, \ldots, x_k= g_k, y_k=h_k $ precisely as in
Theorem \ref{thm:eq}. Then the formula
$$
U(g_1,h_1, \ldots, g_k,h_k) =1 \wedge V(g_1,h_1, \ldots, g_k,h_k) \neq 1
$$
holds in $\GG$. In particular, $U(g_1,h_1, \ldots, g_k,h_k) =1$ in $\GG$. It follows from Corollary \ref{co:merz} that there are words $q_1(x_1) \in \GG[x_1], \ldots, q_k(x_1, \ldots, x_k) \in \GG[x_1,\ldots,x_k]$ such that
$$
\GG[X] \models U(x_1,q_1(x_1, \ldots,x_k), \ldots, x_k,q_k(x_1, \ldots,x_k)) =1.
$$
Moreover, it follows from the construction that $h_1 = q_1(g_1), \ldots, h_k = q_k(g_1, \ldots,g_k)$. We claim that
$$
\GG[X] \models V(x_1,q_1(x_1, \ldots,x_k), \ldots, x_k,q_k(x_1, \ldots,x_k)) \neq 1.
$$
Indeed, if
$$
V(x_1,q_1(x_1, \ldots,x_k), \ldots, x_k,q_k(x_1, \ldots,x_k)) = 1
$$
in $\GG[X]$, then its image in $\GG$ under any specialization $X \rightarrow \GG$ is also trivial, but this is not the case for the specialization  $x_1 \rightarrow g_1, \ldots, x_k \rightarrow g_k$ --- a contradiction. This proves the theorem for  sentences $\phi$ of the form (\ref{eq:sent}). A similar argument works for formulas of the type
$$
\phi = \forall x_1\exists y_1\ldots\forall x_k\exists y_k \bigvee_{i = 1}^{n}(U_i(x_1,y_1, \ldots, x_k,y_k) =1  \wedge V_i(x_1,y_1, \ldots, x_k,y_k) \neq 1),
$$
which is, actually, the general case by Corollary \ref{co:atomic}.
\end{proof}

\section*{Appendix: Positive theory of the direct product of groups}

In this section we prove that if $G=H_1\times \cdots \times H_k$, then the positive theory of $G$ in the language $L_G$ (and in $L$) is decidable if the positive theories of $H_1, \dots, H_k$ are decidable. Perhaps, this result is known, nevertheless, we were not able to find a reference.

The following theorem is due to Feferman and Vaught, see \cite{FV}.

\addtocounter{section}{1}

\begin{theorem} \label{thm:elthdir}
Let $G=H_1\times \cdots \times H_k$. Then the elementary theory of $G$ in the language $L_G$ is decidable, provided that the elementary theory of $H_i$ is decidable in the language $L_{H_i}$, $i=1,\dots, k$.
\end{theorem}
\begin{proof}
Without loss of generality we may assume that $G=A\times B$.

We use induction on the complexity of the formula to prove the following statement.
Given a formula $\varphi(x_1,\dots, x_n)$ in the language $L_G$, one can effectively construct a finite family of formulas $\langle \phi\rangle =\{(\psi_i(y_1,\dots,y_n),\psi_i'(z_1,\dots,z_n))\mid i\in I\}$ such that for all $a_1,\dots, a_n, b_1,\dots, b_n$ we have
$$
A\times B\models \varphi((a_1,b_1), \dots, (a_n,b_n))
$$
if and only if  there exists $i\in I$ such that
$$
A\models \psi_i(a_1,\dots,a_n) \hbox{ and } B\models \psi_i'(b_1,\dots,b_n).
$$
\begin{itemize}
    \item Let $\varphi=(x_i=x_j)$, set $\langle \varphi\rangle=\{(y_i=y_j, z_i=z_j)\}$.
    \item Let $\varphi=(x_i=c)$, where $c\in G$, $c=(c_1, c_2)$, set $\langle \varphi\rangle=\{(y_i=c_1, z_i=c_2)\}$.
    \item Let $\varphi=\varphi_1\vee \varphi_2$ and set $\langle \varphi\rangle=\langle \varphi_1\rangle \cup \langle \varphi_2\rangle$.
    \item Let $\varphi=\neg\varphi_0$ and set
$$
\langle \varphi \rangle =\left\{\left.\left(\bigwedge\limits_{j\in J}\neg \psi_j, \bigwedge\limits_{i\in I\setminus J}\neg \psi_i'\right) \right| J\in \mathcal{P}(I)\right\},
$$
where $\mathcal{P}(I)$ is the power set of $I$ and $\langle \phi_0\rangle =\{(\psi_i,\psi_i')\mid i\in I\}$.
    \item Let $\varphi=\exists x_0 \varphi_0(x_0,x_1,\dots, x_n)$ and set
$$
\langle \varphi \rangle=\{(\exists y_0 \psi_i(y_0,y_1,\dots, y_n),\exists z_0\psi_i' (z_0,z_1,\dots, z_n))
\mid i \in I\},
$$
where $\langle \phi_0\rangle =\{(\psi_i,\psi_i')\mid i\in I\}$.
\end{itemize}
\end{proof}

\begin{cor} \label{cor:pthdir}
Let $G=H_1\times \cdots \times H_k$. Then the positive theory of $G$ in the language $L_G$ {\rm(}in the language $L${\rm)} is decidable, provided that the positive theories of $H_1, \dots, H_k$ are decidable.
\end{cor}
\begin{proof}
In the notation of Theorem \ref{thm:elthdir}, we are left to show that if $\varphi$ is a positive formula in $L_G$ then for all $i\in I$ the formulas $\psi_i$ and $\psi_i'$ are also positive.

By construction of $\langle  \varphi\rangle $ it follows that $\psi_i$ and $\psi_i'$ are positive when
$\varphi= (x_i=x_j)$, $\varphi= (x_i=c)$, $\varphi= P(x_1,\dots, x_n)$, $\varphi=\varphi_1\vee\varphi_2$,  and $\varphi= \exists x_0 \varphi_0(x_0, x_1,\dots, x_n)$. We are left to consider the two following cases: $\varphi=\varphi_1\wedge\varphi_2$ and $\varphi= \forall x_0 \varphi_0(x_0, x_1,\dots, x_n)$.

Let $\varphi= \forall x_0 \varphi_0$. Then $\varphi$ is equivalent to $\neg (\exists x_0\neg\varphi_0(x_0, x_1,\dots, x_n))$. Thus,

\begin{gather} \notag
\begin{split}
\langle \varphi\rangle &= \neg\left\{ \left.\left( \exists x_0\bigwedge\limits_{j\in J} \neg \psi_j,\exists x_0\bigwedge\limits_{i\in I\setminus J} \neg \psi_i' \right)\right| J\in \mathcal{P}(I)\right\}=\\
&\left\{\left.\left(\bigwedge\limits_{J\in J'}\neg\left(\exists x_0\bigwedge\limits_{j\in J}\neg \psi_j\right),\bigwedge\limits_{I\in \mathcal{P} \setminus J'} \neg \left(\exists x_0\bigwedge\limits_{i\in I}\neg\psi_i'\right)\right)\right| J'\in \mathcal{P}(\mathcal{P}(I))\right\}=\\
&
\left\{\left.\left(\bigwedge\limits_{J\in J'}\forall x_0\bigvee\limits_{j\in J} \psi_j,\bigwedge\limits_{I\in \mathcal{P} \setminus J'} \forall x_0\bigvee\limits_{i\in I}\psi_i'\right)\right| J'\in \mathcal{P}(\mathcal{P}(I))\right\}.
\end{split}
\end{gather}

Let now $\varphi=\varphi_1\wedge\varphi_2$ and $\langle \varphi_l \rangle =\{(\psi_{l,i},\psi_{l,i}')\mid i\in I_l\}$, $l=1,2$. Then $\varphi$ is equivalent to $\neg(\neg\varphi_1\vee \neg\varphi_2)$. Thus,

\begin{gather}\notag
\begin{split}
\langle \varphi \rangle&= \neg \left\{\left\{\left.\left(\bigwedge\limits_{j\in J}\neg \psi_{1,j}, \bigwedge\limits_{i\in I_1\setminus J}\neg \psi_{1,i}'\right) \right| J\in \mathcal{P}(I_1)\right\}\bigcup \left\{\left.\left(\bigwedge\limits_{j\in J}\neg \psi_{2,j}, \bigwedge\limits_{i\in I_2\setminus J}\neg \psi_{2,i}'\right) \right| J\in \mathcal{P}(I_2)\right\}\right\}=\\
& \left\{\neg\left\{\left.\left(\bigwedge\limits_{j\in J}\neg \psi_{1,j}, \bigwedge\limits_{i\in I_1\setminus J}\neg \psi_{1,i}'\right) \right| J\in \mathcal{P}(I_1)\right\}\bigcap\neg\left\{\left.\left(\bigwedge\limits_{j\in J}\neg \psi_{2,j}, \bigwedge\limits_{i\in I_2\setminus J}\neg \psi_{2,i}'\right) \right| J\in \mathcal{P}(I_2)\right\}\right\}=\\
&\left\{\left.\left(\bigwedge\limits_{J\in J'}\, \neg \left(\bigwedge\limits_{j\in J}\neg \psi_{1,j}\right),
\bigwedge\limits_{I\in \mathcal{P}(I_1)\setminus J}\,\neg \left(\bigwedge\limits_{i\in I_1\setminus J}\neg \psi_{1,i}'\right)\right) \right| J'\in \mathcal{P}(\mathcal{P}(I_1))\right\}\bigcap\\
&\left\{\left.\left(\bigwedge\limits_{J\in J'}\,\neg \left(\bigwedge\limits_{j\in J}\neg \psi_{2,j}\right), \bigwedge\limits_{I\in \mathcal{P}(I_2)\setminus J}\,\neg \left(\bigwedge\limits_{i\in I_2\setminus J}\neg \psi_{2,i}'\right)\right) \right| J'\in \mathcal{P}(\mathcal{P}(I_2))\right\}=\\
&\left\{\left.\left(\bigwedge\limits_{J\in J'}\, \bigvee\limits_{j\in J} \psi_{1,j},
\bigwedge\limits_{I\in \mathcal{P}(I_1)\setminus J}\, \bigvee\limits_{i\in I_1\setminus J} \psi_{1,i}'\right) \right| J'\in \mathcal{P}(\mathcal{P}(I_1))\right\}\bigcap\\
&\left\{\left.\left(\bigwedge\limits_{J\in J'}\, \bigvee\limits_{j\in J} \psi_{2,j},
\bigwedge\limits_{I\in \mathcal{P}(I_2)\setminus J}\, \bigvee\limits_{i\in I_2\setminus J} \psi_{2,i}'\right) \right| J'\in \mathcal{P}(\mathcal{P}(I_2))\right\}
\end{split}
\end{gather}

\end{proof}

\textbf{Acknowledgements.} The authors are extremely grateful to Andrew Duncan whose numerous remarks and suggestions helped to significantly improve the exposition. We would like to thank Alexei Miasnikov for useful discussions and remarks, and Vladimir Remeslennikov for his comments on an earlier version of this paper.

\end{document}